\documentclass[12pt]{article}
\usepackage{latexsym, amssymb}

\DeclareMathAlphabet\gothic{U}{euf}{m}{n}

\makeatletter
\def\eqnarray{\stepcounter{equation}\let\@currentlabel=\theequation
\global\@eqnswtrue
\tabskip\@centering\let\\=\@eqncr
$$\halign to \displaywidth\bgroup\hfil\global\@eqcnt\z@
  $\displaystyle\tabskip\z@{##}$&\global\@eqcnt\@ne
  \hfil$\displaystyle{{}##{}}$\hfil
  &\global\@eqcnt\tw@ $\displaystyle{##}$\hfil
  \tabskip\@centering&\llap{##}\tabskip\z@\cr}

\def\endeqnarray{\@@eqncr\egroup
      \global\advance\c@equation\m@ne$$\global\@ignoretrue}

\def\@yeqncr{\@ifnextchar [{\@xeqncr}{\@xeqncr[5pt]}}
\makeatother

\begin{document}
\bibliographystyle{tom}

\newtheorem{lemma}{Lemma}[section]
\newtheorem{thm}[lemma]{Theorem}
\newtheorem{cor}[lemma]{Corollary}
\newtheorem{voorb}[lemma]{Example}
\newtheorem{rem}[lemma]{Remark}
\newtheorem{prop}[lemma]{Proposition}
\newtheorem{stat}[lemma]{{\hspace{-5pt}}}
\newtheorem{obs}[lemma]{Observation}
\newtheorem{defin}[lemma]{Definition}

\newenvironment{remarkn}{\begin{rem} \rm}{\end{rem}}
\newenvironment{exam}{\begin{voorb} \rm}{\end{voorb}}
\newenvironment{defn}{\begin{defin} \rm}{\end{defin}}
\newenvironment{obsn}{\begin{obs} \rm}{\end{obs}}

\newenvironment{emphit}{\begin{itemize} }{\end{itemize}}

\newcommand{\gota}{\gothic{a}}
\newcommand{\gotb}{\gothic{b}}
\newcommand{\gotc}{\gothic{c}}
\newcommand{\gote}{\gothic{e}}
\newcommand{\gotf}{\gothic{f}}
\newcommand{\gotg}{\gothic{g}}
\newcommand{\gothh}{\gothic{h}}
\newcommand{\gotk}{\gothic{k}}
\newcommand{\gotm}{\gothic{m}}
\newcommand{\gotn}{\gothic{n}}
\newcommand{\gotp}{\gothic{p}}
\newcommand{\gotq}{\gothic{q}}
\newcommand{\gotr}{\gothic{r}}
\newcommand{\gots}{\gothic{s}}
\newcommand{\gotu}{\gothic{u}}
\newcommand{\gotv}{\gothic{v}}
\newcommand{\gotw}{\gothic{w}}
\newcommand{\gotz}{\gothic{z}}
\newcommand{\gotA}{\gothic{A}}
\newcommand{\gotB}{\gothic{B}}
\newcommand{\gotG}{\gothic{G}}
\newcommand{\gotL}{\gothic{L}}
\newcommand{\gotS}{\gothic{S}}
\newcommand{\gotT}{\gothic{T}}

\newcommand{\mn}{\marginpar{\hspace{1cm}*} }
\newcommand{\mnn}{\marginpar{\hspace{1cm}**} }

\newcommand{\mnq}{\marginpar{\hspace{1cm}*???} }
\newcommand{\mnnq}{\marginpar{\hspace{1cm}**???} }

\newcounter{teller}
\renewcommand{\theteller}{\Roman{teller}}
\newenvironment{tabel}{\begin{list}%
{\rm \bf \Roman{teller}.\hfill}{\usecounter{teller} \leftmargin=1.1cm
\labelwidth=1.1cm \labelsep=0cm \parsep=0cm}
                      }{\end{list}}

\newcounter{tellerr}
\renewcommand{\thetellerr}{(\roman{tellerr})}
\newenvironment{subtabel}{\begin{list}%
{\rm  (\roman{tellerr})\hfill}{\usecounter{tellerr} \leftmargin=1.1cm
\labelwidth=1.1cm \labelsep=0cm \parsep=0cm}
                         }{\end{list}}
\newenvironment{ssubtabel}{\begin{list}%
{\rm  (\roman{tellerr})\hfill}{\usecounter{tellerr} \leftmargin=1.1cm
\labelwidth=1.1cm \labelsep=0cm \parsep=0cm \topsep=1.5mm}
                         }{\end{list}}

\newcommand{\Ni}{{\bf N}}
\newcommand{\Ri}{{\bf R}}
\newcommand{\Si}{{\bf S}}
\newcommand{\Ci}{{\bf C}}
\newcommand{\Ti}{{\bf T}}
\newcommand{\Zi}{{\bf Z}}
\newcommand{\Fi}{{\bf F}}

\newcommand{\proof}{\mbox{\bf Proof} \hspace{5pt}} 
\newcommand{\remark}{\mbox{\bf Remark} \hspace{5pt}}
\newcommand{\ruimte}{\vskip10.0pt plus 4.0pt minus 6.0pt}

\newcommand{\simh}{{\stackrel{{\rm cap}}{\sim}}}
\newcommand{\ad}{{\mathop{\rm ad}}}
\newcommand{\Ad}{{\mathop{\rm Ad}}}
\newcommand{\Aut}{\mathop{\rm Aut}}
\newcommand{\arccot}{\mathop{\rm arccot}}
\newcommand{\capp}{{\mathop{\rm cap}}}
\newcommand{\rcapp}{{\mathop{\rm rcap}}}
\newcommand{\diam}{\mathop{\rm diam}}
\newcommand{\divv}{\mathop{\rm div}}
\newcommand{\codim}{\mathop{\rm codim}}
\newcommand{\RRe}{\mathop{\rm Re}}
\newcommand{\IIm}{\mathop{\rm Im}}
\newcommand{\Tr}{{\mathop{\rm Tr}}}
\newcommand{\Vol}{{\mathop{\rm Vol}}}
\newcommand{\card}{{\mathop{\rm card}}}
\newcommand{\supp}{\mathop{\rm supp}}
\newcommand{\sgn}{\mathop{\rm sgn}}
\newcommand{\essinf}{\mathop{\rm ess\,inf}}
\newcommand{\esssup}{\mathop{\rm ess\,sup}}
\newcommand{\Int}{\mathop{\rm Int}}
\newcommand{\Leibniz}{\mathop{\rm Leibniz}}
\newcommand{\lcm}{\mathop{\rm lcm}}
\newcommand{\loc}{{\rm loc}}

\newcommand{\mod}{\mathop{\rm mod}}
\newcommand{\spann}{\mathop{\rm span}}
\newcommand{\one}{1\hspace{-4.5pt}1}

\newcommand{\DWR}{}

\hyphenation{groups}
\hyphenation{unitary}

\newcommand{\tfrac}[2]{{\textstyle \frac{#1}{#2}}}

\newcommand{\cb}{{\cal B}}
\newcommand{\cc}{{\cal C}}
\newcommand{\cd}{{\cal D}}
\newcommand{\ce}{{\cal E}}
\newcommand{\cf}{{\cal F}}
\newcommand{\ch}{{\cal H}}
\newcommand{\ci}{{\cal I}}
\newcommand{\ck}{{\cal K}}
\newcommand{\cl}{{\cal L}}
\newcommand{\cm}{{\cal M}}
\newcommand{\cn}{{\cal N}}
\newcommand{\co}{{\cal O}}
\newcommand{\cs}{{\cal S}}
\newcommand{\ct}{{\cal T}}
\newcommand{\cx}{{\cal X}}
\newcommand{\cy}{{\cal Y}}
\newcommand{\cz}{{\cal Z}}

\newcommand{\wtozp}{W^{1,2}\raisebox{10pt}[0pt][0pt]{\makebox[0pt]{\hspace{-34pt}$\scriptstyle\circ$}}}
\newlength{\hightcharacter}
\newlength{\widthcharacter}
\newcommand{\covsup}[1]{\settowidth{\widthcharacter}{$#1$}\addtolength{\widthcharacter}{-0.15em}\settoheight{\hightcharacter}{$#1$}\addtolength{\hightcharacter}{0.1ex}#1\raisebox{\hightcharacter}[0pt][0pt]{\makebox[0pt]{\hspace{-\widthcharacter}$\scriptstyle\circ$}}}
\newcommand{\cov}[1]{\settowidth{\widthcharacter}{$#1$}\addtolength{\widthcharacter}{-0.15em}\settoheight{\hightcharacter}{$#1$}\addtolength{\hightcharacter}{0.1ex}#1\raisebox{\hightcharacter}{\makebox[0pt]{\hspace{-\widthcharacter}$\scriptstyle\circ$}}}
\newcommand{\scov}[1]{\settowidth{\widthcharacter}{$#1$}\addtolength{\widthcharacter}{-0.15em}\settoheight{\hightcharacter}{$#1$}\addtolength{\hightcharacter}{0.1ex}#1\raisebox{0.7\hightcharacter}{\makebox[0pt]{\hspace{-\widthcharacter}$\scriptstyle\circ$}}}

 \thispagestyle{empty}

\vspace*{1cm}
\begin{center}
{\Large{\bf The limitations of the Poincar{\'e} inequality}}\\[2mm] 
\large  Derek W. Robinson$^1$ and Adam Sikora$^2$\\[2mm]

\normalsize{May 2013}
\end{center}

\vspace{5mm}

\begin{center}
{\bf Abstract}
\end{center}

\begin{list}{}{\leftmargin=1.8cm \rightmargin=1.8cm \listparindent=10mm 
   \parsep=0pt}
   \item
We examine the validity of the Poincar\'e inequality 
for degenerate, second-order,  elliptic  operators $H$ in divergence form on $L_2(\Ri^{n}\times\Ri^{m})$. 
We assume  the coefficients are real symmetric and $a_1H_\delta\geq H\geq a_2H_\delta$ for some
$a_1,a_2>0$ where $H_\delta$ is a generalized Gru\v{s}in operator,
\[
H_\delta=-\nabla_{x_1}\,|x_1|^{(2\delta_1,2\delta_1')}\,\nabla_{x_1}-|x_1|^{(2\delta_2,2\delta_2')}\,\nabla_{x_2}^2
\;.
\]
Here $x_1\in\Ri^n$, $x_2\in\Ri^m$, $\delta_1,\delta_1'\in[0,1\rangle$,
 $\delta_2,\delta_2'\geq0$
 and $|x_1|^{(2\delta,2\delta')}=|x_1|^{2\delta}$ if $|x_1|\leq 1$ and 
$|x_1|^{(2\delta,2\delta')}=|x_1|^{2\delta'}$ if $|x_1|\geq 1$.

\smallskip

We prove that the Poincar\'e inequality,  formulated in terms of the Riemannian geometry corresponding to $H$, is valid if $n\geq 2$, or if $n=1$ and   $\delta_1\vee\delta_1'\in[0,1/2\rangle$ but it fails if  $n=1$ and  $\delta_1\vee\delta_1'\in[1/2,1\rangle$.
The failure is caused by the leading term. If $\delta_1\in[1/2, 1\rangle$ it is an effect of  the local degeneracy $|x_1|^{2\delta_1}$ but  if  $\delta_1\in[0, 1/2\rangle$ and $\delta_1'\in  [1/2,1\rangle$ it is an effect of the  growth at infinity of 
$|x_1|^{2\delta_1'}$.

If $n=1$ and $\delta_1\in[1/2, 1\rangle$ then the  semigroup $S$ generated by the Friedrichs' extension of $H$ is not ergodic. The subspaces $x_1\geq 0$ and $x_1\leq 0$ are $S$-invariant  and the Poincar\'e inequality is valid on each of these subspaces.
If, however, $n=1$, $\delta_1\in[0, 1/2\rangle$ and $\delta_1'\in  [1/2,1\rangle$ then the semigroup $S$ is ergodic but the Poincar\'e inequality is only valid locally.

\smallskip

Finally we discuss the implication of these results for the kernel of the semigroup $S$.

\vfill
\end{list}

%\vspace{1cm}

\vfill


\bigskip
%\vspace{1cm}
%
%\vspace{1cm}

\noindent
\begin{tabular}{@{}cl@{\hspace{10mm}}cl}
1. & Centre for Mathematics & 
  2. & Department of Mathematics  \\
&\hspace{15mm} and its Applications  & 
  &Macquarie University  \\
& Mathematical Sciences Institute & 
  & Sydney, NSW 2109  \\
& Australian National University& 
  & Australia \\
& Canberra, ACT 0200  & {}
  & \\
& Australia & {}
  & \\
& derek.robinson@anu.edu.au & {}
  &sikora@ics.mq.edu.au \\
\end{tabular}

\newpage
\setcounter{page}{1}

\section{Introduction}\label{S1}

In this paper we continue the analysis of the  class of degenerate  elliptic operators in
divergence form   introduced in \cite{RSi2}. 
The evolution determined by these operators describes diffusion  around and across the surface in $\Ri^d$ on which the corresponding flows are degenerate.
If the degeneracy surface has codimension one then several interesting phenomena can occur depending on the nature of the degeneracy.
In \cite{RSi2} it was established that for sufficiently strong degeneracy  ergodicity can fail; the diffusion can have non-trivial invariant subspaces.
In this situation discontinuous and non-Gaussian behaviour  occurs.
In this paper we will demonstrate that non-Gaussian behaviour can also occur even in the ergodic situation;
the heat kernel  satisfies Gaussian upper bounds but  the matching  lower bounds  are
not necessarily valid.
We will, however,  establish continuity properties and Gaussian bounds for most situations by 
combining the results of \cite{RSi2} with the criteria of 
Grigory'an \cite{Gri4} and  Saloff-Coste  \cite{Sal2, Sal3, Sal4}.
These authors show that Gaussian upper and lower bounds follow from 
two geometric properties, the Poincar\'e inequality and volume doubling.
The crucial feature is that the latter properties are equivalent to the parabolic Harnack inequality of Moser~\cite{Mos}.
Since the volume doubling property was established for the class of  operators we consider by \cite{RSi2}, Corollary~5.2, 
the validity  of lower Gaussian bounds hinges on the Poincar\'e inequality.
The latter property is the main focus of the subsequent analysis.
We demonstrate that the validity of the Poincar\'e inequality is dependent on the order of degeneracy.

The operators we examine  are formally expressed on $\Ri^d$ by
\begin{equation}
H=-\sum^d_{i,j=1}\partial_i\,c_{ij}\,\partial_j
\;,
\label{epi1.1}
\end{equation}
where $\partial_i=\partial/\partial x_i$, the $c_{ij}$ are real-valued measurable functions
and  the  coefficient matrix $C=(c_{ij})$ is symmetric and  positive-definite almost-everywhere.
We assume that $d=n+m$ and adopt the notation $x=(x_1,x_2)$ with $x_1\in \Ri^n$ and $x_2\in\Ri^m$.
Further we assume that 
 $C\sim C_\delta$ where $C_\delta$  is  a block diagonal matrix, 
$C_\delta(x)=|x_1|^{(2\delta_1, 2\delta'_1)}\, I_n+|x_1|^{(2\delta_2, 2\delta'_2)}\,I_m$,
with  $I_n$ and $I_m$ the identity matrices on $\Ri^n$ and $\Ri^m$, respectively.
The indices $\delta_1,\delta_2,\delta'_1,\delta'_2$ are all non-negative,  $\delta_1,\delta_1'<1$
and we use  the notation, introduced in \cite{RSi2}, that $a^{(\alpha,\alpha')}= a^\alpha$ if $ a \leq 1 $
and  $a^{(\alpha,\alpha')}=  a^{\alpha'}$ if $   a \geq 1  $.
Moreover, the equivalence relation $f\sim g$ for  two functions $f$, $g$ with values in a real ordered space indicates that 
there are $a$, $a'>0$ such that $a\,f\leq g\leq a'\,f$. 

The operators are defined more precisely through the  quadratic forms $h$ and $h_\delta$  given by $D(h)=C_c^\infty(\Ri^d)=D(h_\delta)$,
\[
h(\varphi)=\sum^d_{i,j=1}\int_{\Ri^{n+m}} dx\,c_{ij}(x)(\partial_i\varphi)(x)\,(\partial_j\varphi)(x)
\]
and 
\[
h_\delta(\varphi)=\int_{\Ri^d}  dx\,|x_1|^{(2\delta_1,2\delta_1')}((\nabla_{x_1}\varphi)(x))^2+
\int_{\Ri^d}  dx\,|x_1|^{(2\delta_2,2\delta_2')}((\nabla_{x_2}\varphi)(x))^2
\;.
\]
It follows that $h_\delta$ is closable (see, for example, \cite{MR} Section~II.2) and since $h\sim h_\delta$, in the sense of ordering of quadratic forms, $h$ is also closable.
But  by standard arguments the closures $\overline h$ and $\overline h_\delta$ 
are Dirichlet forms.
 Then $H$ and $H_\delta$ are defined as the positive self-adjoint operators associated with these Dirichlet  forms.
  Formally $H_\delta$ is given by
\begin{equation}
H_\delta=-\nabla_{x_1}\,|x_1|^{(2\delta_1,2\delta_1')}\,\nabla_{x_1}-|x_1|^{(2\delta_2,2\delta_2')}\,\nabla_{x_2}^{\,2}
\label{epi1.3}
\end{equation}
and  $H_\delta\sim H$ in the sense of ordering of positive self-adjoint operators.
In the analysis of the  degenerate elliptic operator $H$ the comparison operator $H_\delta$ plays a role analogous to that of the Laplacian in the framework of strongly elliptic operators.
The first  problem  is to establish  properties of the  operators $H_\delta$
and  the second problem is to show that the equivalence relation $H_\delta\sim H$  implies that these properties are shared by $H$.

\smallskip

The Poincar\'e inequality is formulated in terms of the Riemannian geometry defined by the metric $C^{-1}$.
The Riemannian distance  $d(\,\cdot\,;\,\cdot\,)$ can be defined in several equivalent ways.
In particular
\begin{equation}
d(x\,;y)=\sup\{\psi(x)-\psi(y): \psi\in C^1(\Ri^d)\,,\,\Gamma(\psi)\leq1\}
\label{epi1.4}
\end{equation}
for all $x,y\in\Ri^d$
where $\Gamma(\psi)=\sum^d_{i,j=1}c_{ij}(\partial_i\psi)(\partial_j\psi)$ denotes the
 {\it carr\'e du champ} associated with $H$.
 Then the Riemannian ball $B(x\,;r)$ centred at $x\in \Ri^d$ with radius $r>0$ is defined by
 $B(x\,;r)=\{y\in\Ri^d: d(x\,;y)<r\}$.
 The volume (Lebesgue measure) of the ball is denoted by $|B(x\,;r)|$.
 In addition if $n=1$  we set $\Omega_+=\{(x_1,x_2)\in \Ri\times\Ri^m: x_1\geq 0\}$,  $\Omega_-=\{(x_1,x_2)\in \Ri\times\Ri^m: x_1\leq 0\}$  and then  define `balls' $B_\pm(x\,;r)$ by $B_\pm(x\,;r)=B(x\,;r)\cap\Omega_\pm$.

\smallskip

Our principal result is the following.

\begin{thm}\label{tpi1.1} $\;${\bf I.}$\;$ If $n\geq 2$, or if $n=1$ and 
 $\delta_1\vee\delta_1'\in[0,1/2\rangle$, then
then there exist $\lambda>0$ and $\kappa\in\langle0,1]$ such that 
\begin{equation}
\int_{B(x;r)}dy\,\Gamma(\varphi)(y)\geq \lambda\,r^{-2}\int_{B(x;\kappa r)}dy\left(\varphi(y)-\langle\varphi\rangle_{B}\right)^2
\label{epi1.5}
\end{equation}
for all $x\in \Ri^{n+m}$, $r>0$ and  $\varphi\in C^1(\Ri^{n+m})$ where $\langle\varphi\rangle_{B}=|B(x\,;\kappa r)|^{-1}\int_{B(x;\kappa r)}dy\,\varphi(y)$.

\smallskip

\noindent{\bf II.}$\;$ If $n=1$ and  $\delta_1\vee\delta_1'\in[1/2,1\rangle$ then 
the uniform Poincar\'e inequality
$(\ref{epi1.5})$ fails.

\smallskip

\noindent{\bf III.}$\;$ If $n=1$ and 
 $\delta_1\in[1/2, 1\rangle$ then there then there exist $\lambda>0$ and $\kappa\in\langle0,1]$ such that 
\begin{equation}
\int_{B_\pm(x;r)}dy\,\Gamma(\varphi)(y)\geq \lambda\,r^{-2}\int_{B_\pm(x;\kappa r)}dy \left(\varphi(y)-\langle\varphi\rangle_{B_\pm}\right)^2
\label{epi1.6}
\end{equation}
for all $x\in \Omega_\pm$, $r>0$ and  $\varphi\in C^1(\Ri^{1+m})$.
\end{thm}

Theorem~\ref{tpi1.1}  will  be established in Section~\ref{S3}.
The proof is based on the straightforward observation in Section~\ref{S2} that the Poincar\'e inequality~(\ref{epi1.5}) for $H$ is equivalent to the analogous inequality for  $H_\delta$.
This  allows  exploitation of the  characterization derived in \cite{RSi2}, Section~5, of the Riemannian geometry defined by the metric $C_\delta^{-1}$.
Note that (\ref{epi1.5}) and (\ref{epi1.6}) are strong forms of the usual Poincar\'e inequality insofar they are uniform for balls of all position and all  radii but they are weak forms insofar they involve a large ball $B(x\,;r)$ on the left hand side but a small ball $B(x\,;\kappa \,r)$  on the right hand side.
It follows, however,  from the work of Jerison \cite{Jer} (see also \cite{Lu2})
that the weak form together with the volume doubling property implies  the stronger version with $\kappa=1$.

The failure of the Poincar\'e inequality in Case~II has  different origins
in each of  the cases $\delta_1\in [1/2, 1\rangle$ and $\delta_1'\in [1/2, 1\rangle$.
In the first  situation it is caused by the local degeneracy $|x_1|^{2\delta_1}$ of the leading term of $H_\delta$. 
Although the  inequality fails on the whole space $\Ri^{n+m}$ it nevertheless  holds on the subspaces $\Omega_\pm$
by the third statement of the theorem.
If, however, $\delta_1\in [0,1/2\rangle$ and $\delta_1'\in [1/2, 1\rangle$ then the inequality fails because of the growth at infinity of the coefficient $|x_1|^{2\delta_1'}$ in the leading term.
Nevertheless the inequality holds for all $R>0$ and all $B(x\,;r)$ with $r\leq R$  but it does not hold uniformly for all $r$.
It should be emphasized that all these conclusions are independent of the degeneracy parameters $\delta_2,\delta_2'$.

The  Poincar\'e inequality is of relevance for the properties of  the heat kernel corresponding to $H$ because in combination with the volume doubling property it implies both H\"older continuity and upper and lower Gaussian bounds.
This observation follows from the work of Grigor'yan \cite{Gri4} and Saloff-Coste \cite{Sal2, Sal3, Sal4} which extends and simplifies earlier arguments of Moser \cite{Mos, Mos1}.
Since $H$ is defined by the Dirichlet form $\overline h$ it  generates a self-adjoint submarkovian semigroup $S$ on $L_2(\Ri^n\times\Ri^m)$
which, by Proposition~3.1 of \cite{RSi2},  is  bounded as an operator from $L_1(\Ri^n\times\Ri^m)$ to $L_\infty(\Ri^n\times\Ri^m)$.
Therefore $S$ is  determined by a positive bounded integral kernel $K$.
The first statement of Theorem~\ref{tpi1.1} then leads to  the following extension of the results obtained in \cite{RSi2}.

\begin{thm}\label{tpi1.2} $\;$Assume $n\geq 2$, or  $n=1$ and 
 $\delta_1\vee\delta_1'\in[0,1/2\rangle$. 
 Then the semigroup kernel $x,y\in\Ri^{n+m}\mapsto K_t(x\,;y)$ is jointly H\"older continuous 
 and there exist $a,\omega, a',\omega'>0$ such that 
 \begin{eqnarray*}
a\,|B(x\,;t^{1/2})|^{-1}\,e^{-\omega\,d(x;y)^2/t}\leq
K_t(x\,;y)\leq a'\,|B(x\,;t^{1/2})|^{-1}\,e^{-\omega'd(x;y)^2/t}
\end{eqnarray*}
for all $x,y\in\Ri^{n+m}$ and $t>0$.
 \end{thm}

If $n=1$ and  $\delta_1\in[1/2, 1\rangle$ then it follows from Remark~6.9 of \cite{RSi2} that $S_t$ leaves the subspaces
$L_2(\Omega_\pm)$ invariant, i.e.\ the semigroup is not ergodic.
Consequently the kernel $K$ is discontinuous and $K_t(x\,;y)=0$ if $x_1>0$ and $y_1<0$.
Nevertheless the submarkovian semigroups $S^{(\pm)}$ obtained by restricting $S$ to $L_2(\Omega_\pm)$ have bounded,
continuous, integral kernels $K^{(\pm)}$ which satisfy similar Gaussian bounds.

\begin{thm}\label{tpi1.3} $\;$Assume $n=1$ and 
 $\delta_1\in[1/2, 1\rangle$.
   Then the kernels $x,y\in\Omega_\pm\mapsto K_t^{(\pm)}(x\,;y)$ are jointly H\"older continuous 
 and there exist $a,\omega, a',\omega'>0$ such that 
  \begin{eqnarray*}
a\,|B_\pm(x\,;t^{1/2})|^{-1}\,e^{-\omega\,d(x;y)^2/t}\leq 
K_t^{(\pm)}(x\,;y)
\leq a'\,|B_\pm(x\,;t^{1/2})|^{-1}\,e^{-\omega'd(x;y)^2/t}
\end{eqnarray*}
 for all $x,y\in\Omega_\pm$ and $t>0$.
 \end{thm}

The foregoing operators  are related to two classes of degenerate operators which have previously received wide attention.
First if   $\delta_1=\delta_1'=0$ then $H$ is equivalent to  an  operator
\[
H_G=-\nabla_{x_1}^2-c(x_1)\nabla_{x_2}^2
\;.
\]
Operators of this form, with $c(x_1)=x_1^{2k}$ and $k\in \Ni$, were  introduced by  Gru\v{s}in \cite{Gru}.
They are  subelliptic operators of  H\"ormander type  \cite{Hor1} and clearly fall within the class 
we consider.
The Poincar\'e inequality (\ref{epi1.5}) was established for  operators of the form $H_G$ by  Franchi, Guti\`errez and Wheeden   \cite{FGW2}.
These authors considered a wide class of coefficients $c$,
including    $c(x_1)\sim |x_1|^{(2\delta_2,2\delta_2')}$ with $\delta_2,\delta_2'\geq0$, but their methods are completely different to the arguments we use for the operators $H$ and $H_\delta$.
The distinctive feature of the latter operators is  the presence of the  coefficient $|x_1|^{(2\delta_1,2\delta_1')}$
in the leading term.
If $n=1$ these coefficients have a dramatic influence as evinced by Theorems~\ref{tpi1.1} and \ref{tpi1.3}.
 The resulting effects can be partially understood through the associated diffusion process.
The operator $H_\delta$ has   a degeneracy $|x_1|^{(\delta_2,\delta_2')}$ in the components of the underlying flow tangential to the hyperplane $x_1=0$ and  an additional degeneracy $|x_1|^{(\delta_1,\delta_1')}$  in the normal component of the  corresponding flow.
If  $\delta_1\in[1/2,1\rangle$ the normal component of the flow, which is not present in the Gru\v{s}in class, leads to an evolution  which is non-ergodic \cite{RSi2}. 
The corresponding diffusion separates into two distinct components $\Omega_\pm$ and the Poincar\'e inequality on $\Ri^{n+m}$ fails.
The somewhat surprising conclusion is that (\ref{epi1.5})  also fails for $n=1$, $\delta_1\in[0,1/2\rangle$ and  $\delta_1'\in[1/2,1\rangle$  although the diffusion is ergodic.
There is, however,  an approximate failure of ergodicity.
For example, if the one-dimensional diffusion process determined by $-d_x\,(1\vee|x|)\,d_x$
 begins at the right (left) of the origin then with large probability it diffuses to infinity on the right (left).  
Therefore the two  half-lines, $x\geq 0$ and $x\leq 0$ are approximately invariant.
This behaviour is analogous to diffusion on manifolds with ends \cite{ChF} \cite{BCF} \cite{Dav15} \cite{GSC}
and  leads to more complicated lower bounds.
This will be discussed in Section~\ref{S5}.

A second  class  of degenerate operators considered by Trudinger \cite{Tru2} and later by Fabes, Kenig and Serapioni
\cite{FKS} expresses the degeneracy in terms of the 
largest and smallest eigenvalues $\mu_M$, $\mu_m$ of the coefficient matrix.
Typically  Poincar\'e and Harnack inequalities, H\"older continuity etc. follow from local integrability of  $\mu_M$ and $\mu_m^{-1}$.
 These conditions place direct restraints on the order of the local degeneracy, e.g.\ 
 for the operators $H$ under consideration the local integrability of   $\mu_m^{-1}$
 would require that $\delta_1\vee\delta_2<n/2$, and this limits the analysis to weakly degenerate operators.
This type of condition not only rules  out operators such as $H_\delta$ with $n=1$ and $\delta_1\in[1/2,1\rangle$ but also rules out simple examples such as the Heisenberg sublaplacian $H_{\rm Heis}=-\partial_1^2-(\partial_2+x_1\,\partial_3)^2$ on $L_2(\Ri^3)$ for which $\mu_m$ is identically zero.

\section{Preliminaries}\label{S2}

In this section we make some preliminary observations which simplify the subsequent discussion of the Poincar\'e inequality (\ref{epi1.5}).
We begin by recalling some standard consequences of  equivalence properties and then  we derive some approximate scaling estimates.

First  we note  that the Poincar\'e inequality for $H$ is equivalent to the Poincar\'e inequality for the comparison operator $H_\delta$.
This equivalence follows by first remarking that both integrals in the  Poincar\'e inequality (\ref{epi1.5}) are  monotonic functions of the radius of the ball $B=B(x\,;r)$.
This is evident for the left hand integral $\int_B\Gamma(\varphi)$ but it is also true for the right hand integral since \begin{equation}
\int_{B(x; r)}dy\,(\varphi(y)-\langle\varphi\rangle_{B})^2=
\inf_{M\in \Ri}\int_{B(x; r)}dy\,(\varphi(y)-M)^2
\;.
\label{epi2.0}
\end{equation}
This identification is a direct consequence of the  identity
\[
\int_{B}dy\,(\varphi(y)-M)^2=\int_{B}dy\,(\varphi(y)-\langle\varphi\rangle_B)^2+|B|\,(M-\langle\varphi\rangle_B)^2
\;.
\]

Secondly, since  $C\sim C_\delta$ one has 
$\Gamma(\varphi)\sim\Gamma_\delta(\varphi)$
where $\Gamma_\delta$ is the {\it carr\'e du champ} associated with $H_\delta$. 
In particular if $a\,C\leq C_\delta\leq b\,C$ then  $a\,\Gamma(\varphi)\leq \Gamma_\delta(\varphi)\leq b\,\Gamma(\varphi)$.
Consequently, the Riemannian distance $d(\,\cdot\,;\,\cdot\,)$ defined by $\Gamma$ is equivalent to the distance  $d_\delta(\,\cdot\,;\,\cdot\,)$ defined by $\Gamma_\delta$.
Specifically, $b^{-1/2}d(x\,;y)\leq d_\delta(x\,;y)\leq a^{-1/2}d(x\,;y)$ for all $x,y\in \Ri^{n+m}$.
Therefore the corresponding balls $B,B_\delta$ satisfy
\[
B_\delta(x\,;b^{-1/2}r)\subseteq B(x\,;r)\subseteq B_\delta(x\,;a^{-1/2}r)
\]
for all $x\in \Ri^{n+m}$ and $r>0$.
The equivalence of the Poincar\'e inequalities for $H$ and $H_\delta$
 now follows from combination of these remarks.
For example, if (\ref{epi1.5}) is valid then
\begin{eqnarray*}
\int_{B_\delta(x;r)}\Gamma_\delta(\varphi)&\geq &a\int_{B(x;a^{1/2}r)} \Gamma(\varphi)\\[5pt]
&\geq& a\,\lambda\,r^{-2}\inf_{M\in\Ri}\int_{B(x;\kappa a^{1/2}r)}(\varphi(y)-M)^2\\[5pt]
&\geq&a\, \lambda\,r^{-2}\inf_{M\in\Ri}\int_{B_\delta(x;\kappa (a/b)^{1/2}r)}(\varphi(y)-M)^2
\end{eqnarray*}
i.e.\ the analogous inequality is valid for the operator $H_\delta$ but with $\lambda$ replaced by $a\,\lambda$ and $\kappa$ replaced by $(a/b)^{1/2}\kappa$.
The converse implication follows by an identical argument.

Thirdly, one can replace the Riemannian distance $d_\delta(\,\cdot\,;\,\cdot\,)$  by any other equivalent distance 
without destroying the equivalence
property of the Poincar\'e inequality.
In Section~\ref{S3}   we  use the distance function $D_\delta(\,\cdot\,;\,\cdot\,)$ introduced in Section~5 of \cite{RSi2}.
This function is not strictly a distance as it does not satisfy the triangle inequality. 
Nevertheless $D_\delta(\,\cdot\,;\,\cdot\,)$  is equivalent to $d_\delta(\,\cdot\,;\,\cdot\,)$, by Proposition~5.1 of \cite{RSi2},
and the triangle inequality is not used in the foregoing discussion of equivalence
of the balls and  the Poincar\'e inequalities.

\smallskip

A key method for analyzing  differential operators 
is scaling.
Clearly if $\delta_i=\delta_i'$ then $|t\,x|^{(2\delta_i,2\delta_i')}=|t\,x|^{2\delta_i}=t^{2\delta_i}|x|^{(2\delta_i,2\delta_i')}$
for all $t\geq0$.
Consequently simple scaling arguments can then be used to analyze $H_\delta$.
If, however, $\delta_i\neq\delta_i'$ one no longer has a scaling identity.
But one does have scaling estimates.

\begin{prop}\label{ppi2.1}
If $s,t>0$ then
\[
2^{-2(\delta\vee\delta')}\,s^{(\delta,\delta')}\,t^{(\delta\vee\delta', \delta\wedge\delta')}\leq (s\,t)^{(\delta,\delta')}
\leq 2^{2(\delta\vee\delta')}\,s^{(\delta,\delta')}\,t^{(\delta\wedge\delta', \delta\vee\delta')}
\]
for all $\delta,\delta'\geq 0$.
\end{prop}
\proof\
The proof can be established in  two steps.
First one argues that  if $s>0$ and $t\in\langle0,1]$ then 
\begin{equation}
 \left.
 \begin{array}{ll}
  2^{-\delta'-\delta}\,s^{(\delta,\delta')} \,t^{\delta'} & {} \\[5pt]
 2^{-2\delta} \,s^{(\delta,\delta')} \,t^\delta & {}
         \end{array} \right\}
\;\;\leq\; (s\,t)^{(\delta,\delta')}\;\leq \;\;
\left\{
\begin{array}{ll}
  \;\;\;2^{2\delta'}\,s^{(\delta,\delta')} \,t^{\delta} & \;\;\;\mbox{if } \;\delta'\geq \delta  \\[5pt]
\;\;\; 2^{\delta'+\delta} \,s^{(\delta,\delta')} \,t^{\delta'} &\;\;\; \mbox{if } \;\delta\geq \delta'  \; 
\label{est2}
         \end{array} \right.
\end{equation}
for all $\delta,\delta'\geq0$.
Secondly,  if 
 $s>0$ and $t\geq 1$ 
then
\begin{equation}
 \left.
 \begin{array}{ll}
  2^{-\delta'-\delta}\,s^{(\delta,\delta')} \,t^{\delta} & {} \\[5pt]
 2^{-\delta'-\delta} \,s^{(\delta,\delta')} \,t^{\delta'} & {}
         \end{array} \right\}
\;\;\leq\; (s\,t)^{(\delta,\delta')}\;\leq \;\;
\left\{
\begin{array}{ll}
  \;\;\;2^{\delta'+\delta}\,s^{(\delta,\delta')} \,t^{\delta'} & \;\;\;\mbox{if } \;\delta'\geq \delta \; \\[5pt]
\;\;\; 2^{\delta'+\delta} \,s^{(\delta,\delta')} \,t^{\delta} &\;\;\; \mbox{if } \;\delta\geq \delta'  \;
\label{est21}
         \end{array} \right.
\end{equation}
for all $\delta,\delta'\geq0$.
The statement of the proposition is an immediate consequence of these bounds.

The proofs of (\ref{est2}) and (\ref{est21}) are elementary. 
First assume $s>0$ and $t\in\langle0,1]$ and consider (\ref{est2}). 
It is evident that 
$s^{(\delta,\delta')}=(s\wedge1)^\delta(s\vee 1)^{\delta'}$.
Consequently one has estimates
$(s/(1+s))^\delta\leq (s\wedge1)^\delta\leq 2^\delta (s/(1+s))^\delta $
and 
$2^{-\delta'}(1+s)^{\delta'}\leq (s\vee 1)^{\delta'}\leq (1+s)^{\delta'}$.
Hence
\begin{equation}
2^{-\delta'}s^\delta(1+s)^{\delta'-\delta}\leq s^{(\delta,\delta')}\leq 2^\delta s^\delta(1+s)^{\delta'-\delta}
\label{epi2.1}
\end{equation}
for all $s>0$ and all $\delta,\delta'\geq0$.
Then replacing  $s$ with $s\,t$ 
one obtains
\[
(s\,t)^{(\delta,\delta')}\geq 2^{-\delta'}(s\,t)^\delta(1+s\,t)^{\delta'-\delta}
\;.
\]
But if $\delta'\geq \delta $  then 
\[
(1+s\,t)^{\delta'-\delta}\geq (t\,(1+s))^{\delta'-\delta}
\]
since  $t\leq 1$.
Therefore 
\begin{eqnarray*}
(s\,t)^{(\delta,\delta')}&\geq & 2^{-\delta'}(s\,t)^\delta (t\,(1+s))^{\delta'-\delta}\\[5pt]
&=&2^{-\delta'}s^\delta (1+s)^{\delta'-\delta}\,t^{\delta'}
\geq 2^{-\delta'-\delta}s^{(\delta,\delta')}\,t^{\delta'}
\end{eqnarray*}
where the last estimate uses (\ref{epi2.1}).
Alternatively if $\delta'\leq \delta $ then
\[
(1+s\,t)^{\delta'-\delta}\geq (1+s)^{\delta'-\delta}(1+t)^{\delta'-\delta}
\;.
\]
Therefore
\begin{eqnarray*}
(s\,t)^{(\delta,\delta')}&\geq & 2^{-\delta'}(s\,t)^\delta (1+s)^{\delta'-\delta}(1+t)^{\delta'-\delta}
\\[5pt]
&=&2^{-\delta'}s^\delta (1+s)^{\delta'-\delta}\,t^{\delta}(1+t)^{\delta'-\delta}
\geq 2^{-2\delta}s^{(\delta,\delta')}\,t^{\delta}
\end{eqnarray*}
where the last estimate uses (\ref{epi2.1}) and  $t\leq 1$.
Combining these conclusions gives the lower bound of (\ref{est2}). 
The upper bound follows analogously using the upper bound in (\ref{epi2.1}).
The proof of (\ref{est21}) is similar.
We omit the details.
\hfill$\Box$

\bigskip

Finally we examine the action of the scaling semigroup
 $t>0\mapsto \sigma_t$ defined  on $\Ri^{n+m}$ by 
\begin{equation}
\sigma_t(x_1,x_2)=(t^{(\alpha, \alpha')}x_1,t^{(\beta, \beta')}x_2)
\label{epi2.2}
\end{equation}
with 
\[
\alpha=(1-\delta_1)^{-1}\;,\;\;\;\;\; \alpha'=(1-\delta_1')^{-1}\;,
\vspace{2mm}
\]
\[
\beta=(1+\delta_2-\delta_1)\,\alpha\;\;\;{\rm and}\;\;\;
\beta'=(1+\delta_2'-\delta_1')\,\alpha'
\;.
\]
The semigroup acts by transposition on $L_2(\Ri^{n+m})$ and we denote the transpose action by
$\tilde\sigma_t$.
Thus
\[
(\tilde\sigma_t\varphi)(x_1,x_2)=\varphi(t^{(\alpha, \alpha')}x_1,t^{(\beta, \beta')}x_2)
\]
for all $\varphi\in L_2(\Ri^{n+m})$. 

For orientation note that if $\delta_1=\delta_1'$ and $\delta_2=\delta_2'$ then $\alpha=\alpha'$, $\beta=\beta'$ and
\begin{eqnarray*}
\Gamma_\delta(\tilde\sigma_t\varphi)(x)&=&t^{2\alpha-2\alpha\delta_1}|\sigma_t(x_1)|^{\delta_1}|(\tilde\sigma_t\partial_{x_1}\varphi)(x)|^2+t^{2\beta-2\alpha\delta_2}|\sigma_t(x_1)|^{\delta_2}|(\tilde\sigma_t\nabla_{x_2}\varphi)(x)|^2\\[5pt]
&=&t^2\,(\tilde\sigma_t\,\Gamma_\delta(\varphi))(x)
\;.
\end{eqnarray*}
This explains the choice of the  scaling parameters.
They  are  chosen to ensure that the {\it carr\'e du champ} scales quadratically.
The situation is more complicated if $\delta_i\neq\delta_i'$ because there is no exact scaling.
Nevertheless the scaling semigroup has an approximate intertwining property.

\begin{prop}\label{ppi2.2}
Let $\widehat\Gamma$ denote the {\it carr\'e du champ} of the operator
$\widehat H_\delta$ with coefficients $|x_1|^{2(\delta_i\vee \delta_i',\delta_i\wedge \delta_i')}$
and $\widetilde \Gamma$ the {\it carr\'e du champ} of the operator
$\widetilde H_\delta$ with coefficients 
 $|x_1|^{2(\delta_i\wedge \delta_i',\delta_i\vee \delta_i')}$.
Then
\[
2^{4\delta_M}\,t^2\,(\tilde\sigma_t\,\widetilde\Gamma_\delta(\varphi))\geq
\Gamma_\delta(\tilde\sigma_t\varphi)
\geq 
2^{-4\delta_M}\,t^2\,(\tilde\sigma_t\,\widehat\Gamma_\delta(\varphi))
\]
where $\delta_M=\max\{\delta_1,\delta_1',\delta_2,\delta_2'\}$.
\end{prop}
\proof\
It follows by the definition of the scaling semigroup that 
\begin{eqnarray}
\Gamma_\delta(\tilde\sigma_t\varphi)(x)
&=&t^{(2\alpha,2\alpha')}\,|x_1|^{(2\delta_1,2\delta_1')}\,|(\tilde\sigma_t\,\partial_{x_1}\varphi)(x)|^2\nonumber\\[5pt]
&&\hspace{2.5cm}{}+
t^{(2\beta,2\beta')}
\,|x_1|^{(2\delta_2,2\delta_2')}|(\tilde\sigma_t\nabla_{x_2}\varphi)(x)|^2\label{eas2.1}\;.
\end{eqnarray}
But  the lower bound in Proposition~\ref{ppi2.1} gives
\begin{eqnarray*}
|x_1|^{(2\delta_i,2\delta_i')}&=&
(t^{(-\alpha, -\alpha')}\,|t^{(\alpha, \alpha')}x_1|)^{(2\delta_i,2\delta_i')}\\[5pt]
&\geq& 2^{-4(\delta_i\vee\delta_i')}\,
\,t^{(-2\alpha\delta_i,-2\alpha'\delta_i')}\,|t^{(\alpha, \alpha')}x_1|^{2(\delta_i\vee\delta_i',\delta_i\wedge\delta_i')}
\end{eqnarray*}
for both  $i=1$ and $i=2$.
Combination of these estimates then gives
\begin{eqnarray*}
\Gamma_\delta(\tilde\sigma_t\varphi)(x)&\geq& 2^{-4(\delta_1\vee\delta_1')}t^{(2\alpha-2\delta_1\alpha,2\alpha'-2\delta_1'\alpha')}\,
|\sigma_t(x_1)|^{2(\delta_1\vee\delta_1',\delta_1\wedge\delta_1')}|(\tilde\sigma_t\,\partial_{x_1}\varphi)(x)|^2\\[5pt]
&&\hspace{2cm}{}+2^{-4(\delta_2\vee\delta_2')}\,
t^{(2\beta-2\alpha\delta_2,2\beta'-2\alpha'\delta_2')}\,|\sigma_t(x_1)|^{2(\delta_2\vee\delta_2',\delta_2\wedge\delta_2')}|(\tilde\sigma_t\nabla_{x_2}\varphi)(x)|^2\\[5pt]
&\geq& 2^{-4(\delta_1\vee\delta_1')}t^2\,
|\sigma_t(x_1)|^{2(\delta_1\vee\delta_1',\delta_1\wedge\delta_1')}|(\tilde\sigma_t\,\partial_{x_1}\varphi)(x)|^2\\[5pt]
&&\hspace{2cm}{}+2^{-4(\delta_2\vee\delta_2')}\,
t^2\,|\sigma_t(x_1)|^{2(\delta_2\vee\delta_2',\delta_2\wedge\delta_2')}|(\tilde\sigma_t\nabla_{x_2}\varphi)(x)|^2
\end{eqnarray*}
because $2\alpha-2\delta_1\alpha=2$, $2\beta-2\alpha\delta_2=2$ etc.
Therefore 
\begin{eqnarray*}
\Gamma_\delta(\tilde\sigma_t\varphi)(x)
&\geq& 2^{-4\delta_M}t^{2}\,(\tilde\sigma_t\,\widehat \Gamma_\delta(\varphi))(x)
\;.
\end{eqnarray*}
The upper bound is derived similarly but using the upper bound of
Proposition~\ref{ppi2.1}.
\hfill$\Box$

\section{Poincar\'e inequality}\label{S3}

In this section we prove 
Theorem~\ref{tpi1.1}.
The main onus of the proof consists of establishing the Poincar\'e inequality (\ref{epi1.5}).
The proof of the analogous inequality (\ref{epi1.6})  for $n=1$ and $\delta_1\in [1/2,1\rangle$ is an almost direct consequence of the argument.

It follows  from the discussion of equivalences in Section~\ref{S2} that it suffices to prove the 
Poincar\'e inequality for the operator $H_\delta$. 
Then since $H_\delta$ is invariant under translations in the~{$x_2$-directions} it is sufficient to consider balls with centres 
$(x_1,0)$.
The proof  will be broken down into three distinct cases.
First we examine  balls centred at the origin $(0,0)$ and secondly
 balls that do not contain the origin. 
 Finally we deduce the result for general balls from the two special cases.

\bigskip

\noindent{\bf Case I--Balls centred at the origin.}
This case is handled in three  steps. 
First we  derive the Poincar\'e inequality for a unit cube centred at the origin.
Secondly we extend the result to parallelepipeds obtained by scaling the cube with the semigroup
of scale transformations introduced in  Section~\ref{S2}.
Finally we establish embedding properties involving the balls and parallelepipeds which allow the deduction of the desired inequality for balls.

\begin{prop}\label{ppi3.1}
Assume $n\geq2$ or $n=1$ and $\delta_1\in[0,1/2\rangle$.
Then there is a $\lambda>0$ such that 
\begin{equation}
\int_{[-1,1]^{n+m}}dx\,\Gamma_\delta(\varphi)(x)\geq\lambda\int_{[-1,1]^{n+m}}dx\,(\varphi(x)-\langle\varphi\rangle)^2
\label{epi3.1}
\end{equation}
for all $\varphi\in C_c^1(\Ri^{n+m})$ where $\langle\varphi\rangle=2^{-(n+m)}\int_{[-1,1]^{n+m}}dx\,\varphi(x)$.
\end{prop}
\proof\
First  since the proposition only involves $x$ with $|x_1|\leq 1$ one can assume the coefficients of $H_\delta$ are given by $|x_1|^{2\delta_i}$.
Thus we may assume that 
\[
\Gamma_\delta(\varphi)(x)=|x_1|^{2\delta_1}\,((\nabla_{x_1}\varphi(x))^2+|x_1|^{2\delta_2}((\nabla_{x_2}\varphi)(x))^2
\;.
\]
Secondly, 
the quadratic form
\[
\varphi\in C_c^\infty(\Ri^{n+m})\mapsto h_\delta(\varphi)=
\int_{[-1,1]^{n+m}}dx\,\Gamma_\delta(\varphi)(x)
\]
is closable (see, for example, \cite{MR} Section~II.2a).
Then by standard arguments the subspace $\cd$ of $C_c^1(\Ri^{n+m})$ consisting of functions whose normal derivative is zero on the boundary of the  parallelepiped $[-1,1]^{n+m}$ is a core of $\overline h_\delta$.
(Formally the closure  $\overline h_\delta$ determines the self-adjoint version of the operator (\ref{epi1.3})
corresponding to Neumann boundary conditions.)
Hence for the first statement of the proposition it suffices to verify (\ref{epi3.1}) on~$\cd$.

Thirdly let $\widetilde \varphi$ denote the Fourier cosine-transform with respect to the $x_2$-variables of $\varphi\in\cd$
and $\widetilde{\varphi}_{x_1}$ the cosine-transform of the gradient $\nabla_{x_1}\varphi$.
Then \begin{eqnarray*}
h_\delta(\varphi)=
\sum_{k\in\Zi^m}\int_{[-1,1]^{n}}dx_1\,\Big(|x_1|^{2\delta_1}(\widetilde{\varphi}_{x_1}(x_1,k))^2
+(\pi k/2)^2\,|x_1|^{2\delta_2}(\widetilde\varphi(x_1,k))^2\Big)
\end{eqnarray*}
and
\[
\int_{[-1,1]^{n+m}}dx\,(\varphi(x)-\langle\varphi\rangle)^2
=\int_{[-1,1]^{n}}dx_1\,\left(\widetilde\varphi(x_1,0)-\langle\widetilde\varphi\rangle_0\right)^2
+\sum_{k\in\Zi^m\backslash\{0\}}\int_{[-1,1]^{n}}dx_1\,(\widetilde\varphi(x_1,k))^2
\]
where 
\[
\langle\varphi\rangle=2^{-n-m}\int_{[-1,1]^{n+m}}dx\,\varphi(x)
=2^{-n}\int_{[-1,1]^{n}}dx_1\,\widetilde\varphi(x_1,0)=\langle\widetilde\varphi\rangle_0\;.
\]
Therefore to establish (\ref{epi3.1}) it suffices to prove that one can choose $\lambda>0$ such that 
\[
\int_{[-1,1]^{n}}dx_1\,|x_1|^{2\delta_1}(\widetilde\varphi_{x_1}(x_1,0))^2
\geq \lambda \int_{[-1,1]^{n}}dx_1\,(\widetilde\varphi(x_1,0)-\langle\widetilde\varphi\rangle_0)^2
\]
and in addition
\[
\int_{[-1,1]^{n}}dx_1\,\Big(|x_1|^{2\delta_1}(\widetilde\varphi_{x_1}(x_1,k))^2
+(\pi k/2)^2\,|x_1|^{2\delta_2}(\widetilde\varphi(x_1,k))^2\Big)
\geq \lambda \int_{[-1,1]^{n}}dx_1\,(\widetilde\varphi(x_1,k))^2
\]
for all $k\in\Zi^m\backslash\{0\}$.

Fourthly, if for $x\in\Ri^n$ and $k\in \Zi^m$  one defines $\psi_k$ by setting $\psi_k(x)=\widetilde\varphi(x,k)$ then $\psi\in C_c^1(\Ri^n)$.
Therefore it now suffices to prove that there is a $\lambda>0$ such that
\begin{equation}
\int_{[-1,1]^{n}}dx\,|x|^{2\delta_1}((\nabla_x\psi)(x))^2\geq \lambda \int_{[-1,1]^{n}}dx\,(\psi(x)-\langle\psi\rangle)^2
\label{epi3.2}
\end{equation}
and, in addition,
\begin{equation}
\tilde h_\delta(\psi)
\geq \lambda \int_{[-1,1]^{n}}dx\,(\psi(x))^2
\label{epi3.3}
\end{equation}
for all $\psi\in C_c^1(\Ri^n)$ where $\langle\psi\rangle=2^{-n}\int_{[-1,1]^{n}}\psi$ 
and $\tilde h_\delta$ denotes the form 
\begin{equation}
\tilde h_\delta(\psi)=\int_{[-1,1]^{n}}dx\,\Big(|x|^{2\delta_1}((\nabla_x\psi)(x))^2
+(\pi/2)^2\,|x|^{2\delta_2}(\psi(x))^2\Big)
\label{epi3.4}
\end{equation}
on $L_2([-1,1]^n)$ with domain given by the restrictions of $C_c^1(\Ri^n)$ to the parallelepiped $[-1,1]^n$.
These two properties are established by the following lemmas.

\begin{lemma}\label{lpi3.1}
If $n\geq2$ or $n=1$ and $\delta_1\in[0,1/2\rangle$ then
 one may choose $\lambda>0$ such that $(\ref{epi3.2})$ is valid.
\end{lemma}
\proof\ Let $\widetilde H$ denote the positive self-adjoint operator associated with the closure of the form
$\psi\in C_c^1(\Ri^n)\mapsto\int_{[-1,1]^{n}}dx\,|x|^{2\delta_1}((\nabla_x\psi)(x))^2$.
It follows by standard arguments that $\widetilde H$ has compact resolvent.
Now zero is an eigenvalue and 
 if $\varphi$ is a corresponding eigenfunction  then 
$\int_{[-1,1]^{n}}dx\,|x|^{2\delta_1}((\nabla_x\varphi)(x))^2=0$. 
Therefore $\varphi=0$ on the complement of the origin.
Thus if $n\geq 2$ one must have $\varphi=0$ and  zero is a simple  eigenvalue.
But if $n=1$ then $\varphi$ is  constant on $[-1,0\rangle$ and on $\langle0,1]$ and the  
zero  eigenvalue has  multiplicity two. 

If $n\geq2$ it follows that there is a $\lambda>0$ such that $\widetilde H\geq \lambda I$ on the orthogonal complement of the constants.
But this is just an alternative formulation of (\ref{epi3.2}).

If $n=1$ the foregoing argument does not work and indeed (\ref{epi3.2}) fails if $\delta_1\in[1/2,1\rangle$.
But if $\delta_1\in[0,1/2\rangle$ and $\psi\in C_c^1(\Ri)$  it 
follows from the Cauchy--Schwarz inequality that
\begin{eqnarray*}
|\psi(x)-\psi(0)|^2&=&\Big|\int^x_0ds\,|s|^{-\delta_1}\Big(|s|^{\delta_1}\psi'(s)\Big)\Big|^2\\[5pt]
&\leq&\int^1_0 ds\,s^{-2\delta_1}\,\int^1_{-1}ds\,|s|^{2\delta_1}(\psi'(s))^2
=(1-2\delta_1)^{-1}\int^1_{-1}ds\,|s|^{2\delta_1}(\psi'(s))^2
\end{eqnarray*}
for all $x\in[-1,1]$.
Therefore, using (\ref{epi2.0}), one has
\begin{eqnarray*}
\int^1_{-1}dx\,|x|^{2\delta_1}(\psi'(x))^2&\geq& ((1-2\delta_1)/2)\int^1_{-1}dx\,(\psi(x)-\psi(0))^2\\[5pt]
&\geq & ((1-2\delta_1)/2)\int^1_{-1}dx\,(\psi(x)-\langle\psi\rangle)^2
\end{eqnarray*}
and  (\ref{epi3.2}) is valid with $\lambda=(1-2\delta_1)/2$.
\hfill$\Box$

\begin{lemma}\label{lpi3.2}
If $n\geq2$ or $n=1$ and $\delta_1\in[0,1/2\rangle$ then
one may choose $\lambda>0$ such that $(\ref{epi3.3})$ is valid.
\end{lemma}
\proof\  
It follows from Lemma~\ref{lpi3.1} that one may choose $\lambda>0$ such that (\ref{epi3.2}) is valid
 for all $\psi\in C_c^1(\Ri^n)$.
Set $R_l=[-l,l\,]^{n}$ with $l\in\langle0,1\rangle$ and let $\|\cdot\|_2$ denote the $L_2(R_1)$-norm.

First assume  $\|\psi-\langle\psi\rangle\|_2^2\geq \|\psi\|_2^2/4$.
Then it follows from  (\ref{epi3.2}) that $\tilde h_\delta(\psi)\geq (\lambda/4)\,\|\psi\|_2^2$.

Secondly, assume $\|\psi-\langle\psi\rangle\|_2^2\leq \|\psi\|_2^2/4$.
Then
\begin{eqnarray*}
 \int_{R_1}dx\,|x|^{2\delta_2}\,|\psi(x)|^2\geq l^{2\delta_2}\int_{R_1\backslash R_l}dx\,|\psi(x)|^2
  =l^{2\delta_2}\Big(\|\psi\|_2^2-\int_{R_l}dx\,|\psi(x)|^2\Big)
 \end{eqnarray*}
 for all $l\in\langle0,1]$.
 But 
\begin{eqnarray*}
\int_{R_l}dx\,|\psi(x)|^2&\leq&2\int_{R_l}|\psi(x)-\langle\psi\rangle|^2+2\,\langle\psi\rangle^2\,|R_l|\\[5pt]
&\leq&2\,\|\psi-\langle\psi\rangle\|_2^2+2\,\|\psi\|_2^2\,(|R_l|/|R_1|)\leq 2\,(1/4+l^n)\,\|\psi\|_2^2
\;.
\end{eqnarray*}
Now combining the last two estimates and setting $l=1/8$ one has
\[
 \int_{R_1}dx\,|x|^{2\delta_2}\,|\psi(x)|^2\geq 
 (2^{-6\delta_2}/4)\,\|\psi\|_2^2
 \;.
 \]
 Then it follows from (\ref{epi3.2})  that  $\tilde h_\delta(\psi)\geq (\pi/2)^2\,(2^{-6\delta_2}/4)\,\|\psi\|_2^2$.

 \smallskip
One concludes that 
  $\tilde h_\delta(\psi)\geq \lambda_0\,\|\psi\|_2^2$ for all $\psi\in C_c^1(\Ri^n)$
 with $\lambda_0=(\lambda\wedge  2^{-6\delta_2})/4$.\hfill$\Box$

\bigskip

The 
statement of Proposition~\ref{ppi3.1} now follows from Lemmas~\ref{lpi3.1} and \ref{lpi3.2} 
by the discussion preceding the lemmas.\hfill$\Box$

\bigskip

Proposition~\ref{ppi3.1} establishes the Poincar\'e inequality on the Euclidean cube $C_1=[-1,1]^{n+m}$ and we now extend the result to more general cubes by the scaling transformation (\ref{epi2.2}) introduced in Section~\ref{S2}.
Let $C_t=\sigma_t(C_1)$ for all $t>0$.
 Explicitly 
 \begin{eqnarray}
 C_t &=&\{x\in \Ri^{n+m}: |x_1|_\infty< t^{(\alpha,\alpha')},\; |x_2|_\infty< t^{(\beta,\beta')} \}\nonumber\\[5pt]
 &=&\{x\in \Ri^{n+m}: (|x_1|_\infty)^{(1-\delta_1, 1-\delta_1')}< t,\; (|x_2|_\infty)^{(1-\gamma, 1-\gamma')}< t \}
 \label{epi3.21}
\end{eqnarray} 
where $\gamma=\delta_2(1+\delta_2-\delta_1)^{-1}$, $\gamma'=\delta'_2(1+\delta'_2-\delta'_1)^{-1}$
and $|x_1|_\infty$,  $|x_2|_\infty$ denote the $l_\infty$-norms of $x_1\in\Ri^n$ and $x_2\in \Ri^m$.

Next we 
 apply Proposition~\ref{ppi3.1} to the operator with the coefficients $2(\delta_i\vee\delta_i',\delta_i\wedge\delta_i')$
to reach the following conclusion.

\begin{prop}\label{ppi3.2}
Assume $n\geq 2$ or $n=1$ and $\delta_1\vee\delta_1'\in[0,1/2\rangle$.
Then there is a $\lambda>0$ such that 
\begin{equation}
\int_{C_t}dx\,\Gamma_\delta(\varphi)(x)\geq\lambda\,t^{-2}\int_{C_t}dx\,(\varphi(x)-\langle\varphi\rangle)^2
\label{epi3.22}
\end{equation}
for all $\varphi\in C_c^1(\Ri^{n+m})$  and $t>0$ where $\langle\varphi\rangle=|C_t|^{-1}\int_{C_t}dx\,\varphi(x)$.
\end{prop}
\proof\
First by a change of coordinates $x\in C_t\to y=\sigma_t(x)\in C_1$ one has
\[
\int_{C_t}dx\,\Gamma_\delta(\varphi)(x)=J(t)\int_{C_1}dy\,(\tilde\sigma_t\,\Gamma_\delta(\varphi))(y)
\]
with $J(t)=t^{(n\alpha+m\beta, n\alpha'+m\beta')}$ the Jacobian of the coordinate change.
Secondly, it follows from   the lower bound of  Proposition~\ref{ppi2.2} that 
$(\tilde\sigma_t\,\Gamma_\delta(\varphi))\geq 2^{-4\delta_M}\,t^{-2} \,\widehat \Gamma_\delta(\tilde\sigma_t\varphi)$.
Therefore
\[
J(t)\int_{C_1}dy\,(\tilde\sigma_t\,\Gamma_\delta(\varphi))(y)
\geq 2^{-4\delta_M}\,J(t)\,t^{-2}\, \int_{C_1}dy\,\widehat \Gamma_\delta(\tilde\sigma_t\varphi)(y)
\;.
\]
Thirdly,  $\widehat \Gamma_\delta$ is the {\it carr\'e du champ} of the operator with coefficients 
$|x_1|^{2(\delta_i\vee \delta_i',\delta_i\wedge \delta_i')}$.
Then since $n\geq2$ or $\delta_1\vee\delta_1'\in[0,1/2\rangle$ if $n=1$
 one can apply the Poincar\'e inequality of Proposition~\ref{ppi3.1} together with the identification (\ref{epi2.0}) to deduce that 
\[
\int_{C_1}dy\,\widehat \Gamma_\delta(\tilde\sigma_t\varphi)(y)
\geq \lambda \inf_{M\in \Ri}\int_{C_1}dy\,\left((\tilde\sigma_t\varphi)(y)-M\right)^2
\;.
\]
Therefore by combination of these observations and another coordinate change one finds
\begin{eqnarray*}
\int_{C_t}dx\,\Gamma_\delta(\varphi)(x)
&\geq &
 \lambda_\delta\,t^{-2}\,J(t)\,\inf_{M\in \Ri}\int_{C_1}dx\,\left((\tilde\sigma_t\varphi)(x)-M\right)^2
\\[5pt]
&=& \lambda_\delta\,t^{-2}\,\inf_{M\in \Ri}\int_{C_t}dx\,\left(\varphi(x)-M\right)^2
=\lambda_\delta\,t^{-2}\int_{C_t}dx\,(\varphi(x)-\langle\varphi\rangle)^2
 \end{eqnarray*}
 for all $t>0$ where $\lambda_\delta= \lambda\,2^{-4\delta_M}$.
 \hfill$\Box$
 
 \bigskip
 
At this point we appeal to the discussion given in Section~5 of \cite{RSi2} of the Riemannian geometry defined by 
the metric $C_\delta^{-1}$.
The corresponding Riemannian distance $d_\delta(\,\cdot\,;\,\cdot\,)$ 
is equivalent to the distance given by the function $D_\delta(\,\cdot\,;\,\cdot\,)$ where 
\begin{eqnarray}
D_\delta(x\,;y)&=&|x_1-y_1|/(|x_1|+|y_1|)^{(\delta_1, \delta_1')}\nonumber\\[5pt]
&&\hspace{1cm}{}+|x_2-y_2|\left((|x_1|+|y_1|)^{(\delta_2,\delta_2')}+
(|x_2|+|y_2|)^{(\gamma, \gamma')}\right)^{-1}\label{epi3.41}
\end{eqnarray}
with   $|x_i|$   the $l_2$-norm of $x_i$.
In fact $D_\delta(\,\cdot\,;\,\cdot\,)$ is not strictly a distance since it does not satisfy the triangle inequality but, as mentioned in Section~\ref{S2}, this does not affect the discussion of the Poincar\'e inequality. 
It suffices that  $D_\delta(\,\cdot\,;\,\cdot\,)$ is equivalent
to the Riemannian distance.

Now we extend the Poincar\'e inequality of Proposition~\ref{ppi3.2} from the parallelepipeds $C_t$ to the centred balls $B_\Delta(0\,;r)$ defined by the distance  function $D_\delta(\,\cdot\,;\,\cdot\,)$.
The extension is based on the following two embedding lemmas.

\begin{lemma}\label{lpi3.3}
$\;\;\; C_t\subseteq B_\Delta(0\,;4\,(n+m)\,t)$ for all $t>0$.
\end{lemma}
\proof\
First if $x_1\in \Ri^n$ then $|x_1|\leq n^{1/2}\,|x_1|_\infty\leq n\,|x_1|_\infty$.
Therefore using  Proposition~\ref{ppi2.1} one has $|x_1|^{(1-\delta_1,1-\delta_1')}\leq 4\,n\,(|x_1|_\infty)^{(1-\delta_1,1-\delta_1')}$.
Similarly  $|x_2|^{(1-\gamma,1-\gamma')}\leq 4\,n\,(|x_2|_\infty)^{(1-\gamma,1-\gamma')}$
for  $x_2\in\Ri^m$.

Secondly, it  follows from the  characterization (\ref{epi3.21}) of $C_t$ that
\begin{eqnarray*}
C_t&\subseteq& \{x\in\Ri^{n+m}:4\,n\,(|x_1|_\infty)^{(1-\delta_1, 1-\delta_1')}+4\,m\,(|x_2|_\infty)^{(1-\gamma, 1-\gamma')}< 4\,(n+m)\,t\}\\[5pt]
&\subseteq&
\{x\in\Ri^{n+m}:|x_1|^{(1-\delta_1, 1-\delta_1')}+|x_2|^{(1-\gamma, 1-\gamma')}< 4\,(n+m)\,t\}\\[5pt]
&\subseteq&
\{x\in\Ri^{n+m}:D_\delta(x\,;0)< 4\,(n+m)\,t\}=B_\Delta(0\,;4\,(n+m)\,t)
\end{eqnarray*}
for all $t>0$ where we have used $
|x_2|\left(|x_1|^{(\delta_2,\delta_2')}+
(|x_2|)^{(\gamma, \gamma')}\right)^{-1}
\leq (|x_2|)^{(1-\gamma, 1-\gamma')}
$.
\hfill$\Box$

\begin{lemma}\label{lpi3.4} 
There is a $\kappa\in\langle0, 1]$ such that $B_\Delta(0\,;\kappa\,t)\subseteq C_t$ for all $t>0$.
\end{lemma}
\proof\
If $x\in B_\Delta(0\,;t)$ then $|x_1|^{(1-\delta_1,1-\delta_1')}<t$ and 
\[
|x_2|<t\,\left(|x_1|^{(\delta_2,\delta_2')}+|x_2|^{(\gamma, \gamma')}\right)\;.
\]
Therefore $|x_1|<t^{(\alpha,\alpha')}$ and 
\begin{eqnarray}
|x_2|&<&t^{(1+\alpha\delta_2,1+\alpha'\delta_2')}+t\,|x_2|^{(\gamma, \gamma')}\nonumber \\[5pt]
&=&t^{(\beta,\beta')}+t\,|x_2|^{(\gamma, \gamma')}\label{epi3.5}
\end{eqnarray}
where $\alpha, \alpha',\beta$ and $\beta'$  are the parameters introduced in the definition (\ref{epi2.2}) of the scaling semigroup.

Now we consider  the cases $t\leq 1$ and $t\geq1$ separately.

First, if $t\leq 1$ then $|x_1|\leq t^\alpha$.
Moreover, 
$|x_2|\leq t^{\beta}+t\,|x_2|^{(\gamma, \gamma')}\leq 1+|x_2|^{(\gamma, \gamma')}$.
Then since $\gamma, \gamma'<1$ it follows that  there is an $a>0$ such that $|x_2| \leq a$.
There are two possibilities, $a\leq 1$ or $a>1$.
If $a\leq 1$ then $|x_2|^{(\gamma, \gamma')}=|x_2|^\gamma$ and 
$|x_2|\leq t^{\beta}+t\,|x_2|^{\gamma}$
or, equivalently,
\[
(|x_2|/t^\beta)\leq 1+(|x_2|/t^\beta)^{\gamma}
\]
for all $t\in\langle0,1]$.
Then since $\gamma<1$  one can choose $b>0$ such that $|x_2|\leq b\,t^\beta$
for all $t\leq 1$.
Alternatively if $a>1$ then
\[
|x_2|^{(\gamma, \gamma')}=(a\,(|x_2|/a))^{(\gamma, \gamma')}
\leq 2^{2(\gamma\vee\gamma')}\,a^{\gamma\vee\gamma'}(|x_2|/a)^{\gamma}\;,
\]
by the upper bound of Proposition~\ref{ppi2.1}, and one now has
\[
(|x_2|/t^\beta)\leq 1+a'\,(|x_2|/t^\beta)^{\gamma}
\]
with $a'=2^{2(\gamma\vee\gamma')}\,a^{\gamma\vee\gamma'-\gamma}$.
Therefore one again deduces a bound $|x_2|\leq b\,t^\beta$
for all $t\leq 1$.
Thus if   $\kappa\leq (1\vee b)^{-\beta}$ then $B_\Delta(0\,;\kappa t)\subseteq C_{t}$ for all $t\leq 1$.

Secondly  suppose $t\geq 1$.
Then it follows from (\ref{epi3.5}) that
\begin{eqnarray*}
(|x_2|/t^{\beta'})&\leq &1+t^{1-\beta'}\,|x_2|^{(\gamma,\gamma')}\\[5pt]
&=& 1+t^{\beta'\gamma'}\,(t^{\beta'}(|x_2|/t^{\beta'}))^{(\gamma,\gamma')}
\leq 1+2^{\gamma+\gamma'}\,(|x_2|/t^{\beta'})^{\gamma\vee\gamma'}
\end{eqnarray*}
by another application of the upper bounds of Proposition~\ref{ppi2.1}.
Since $\gamma\vee\gamma'<1$ it follows that there is a $b'>0$ such that $|x_2|\leq b'\,t^{\beta'}$
uniformly for all $t\geq 1$.
But one also has the bound $|x_1|<t^{\alpha'}$ for all $t\geq 1$.
(This is evident if $ |x_1|\leq 1$ but if $|x_1|\geq 1$ then $|x_1|^{1-\delta_1'}\leq t$ and the bound  again follows.)
Therefore one now concludes that if $\kappa\leq (1\vee b')^{-\beta}$ then $B_\Delta(0\,;\kappa' t)\subseteq C_{t}$ for all $t\geq 1$.
The statement of the lemma follows immediately.
\hfill$\Box$

\bigskip

The Poincar\'e inequality now extends to the balls $B_\Delta$.

\begin{prop}\label{ppi3.3}
Assume $n\geq 2$ or $n=1$ and $\delta_1\vee\delta_1'\in[0,1/2\rangle$.
Then there are $\lambda_1>0$ and  $\kappa_1\in\langle0,1]$  such that 
\begin{equation}
\int_{B_\Delta(0;r)}dx\,\Gamma_\delta(\varphi)(x)\geq\lambda_1\,r^{-2}\int_{B_\Delta(0;\kappa_1 r)}dx\,(\varphi(x)-\langle\varphi\rangle)^2
\label{epi3.33}
\end{equation}
for all $\varphi\in C_c^1(\Ri^{n+m})$  and $r>0$ where $\langle\varphi\rangle$ is the average of $\varphi$ over $B_\Delta(0;\kappa_1 r)$.
\end{prop}
\proof\
Set $\hat r=r/(4(n+m))$.
It follows from Proposition~\ref{ppi3.2} together with Lemmas~\ref{lpi3.3} and \ref{lpi3.4} that 
\begin{eqnarray*}
\int_{B_\Delta(0;r)}dx\,\Gamma_\delta(\varphi)(x)&\geq &
\int_{C_{\hat r}}dx\,\Gamma_\delta(\varphi)(x)\\[5pt]
&\geq&\lambda\,{\hat r}^{-2}\inf_{M\in\Ri}\int_{C_{\hat r}}dx\,(\varphi(x)-M)^2\\[5pt]
&\geq&16\,(n+m)^2\,\lambda\,r^{-2}\inf_{M\in\Ri}\int_{B_\Delta(0;\kappa \hat r)}dx\,(\varphi(x)-M)^2\\[5pt]
&=&16\,(n+m)^2\,\lambda\,r^{-2}\int_{B_\Delta(0;\kappa_1 r)}dx\,(\varphi(x)-\langle\varphi\rangle)^2
\end{eqnarray*}
which gives the desired conclusion with $\lambda_1=16(n+m)^2\,\lambda$ and $\kappa_1= \kappa/(4(n+m))$.
\hfill$\Box$

\bigskip

The last proposition  establishes the Poincar\'e inequality for the Riemannian balls $B_\Delta(0\,;r)$ for all $r>0$.

\bigskip

\noindent{\bf Case II--Balls not containing the origin.}

Next we consider balls  $B_\Delta((\xi_1,0)\,;r)$ determined by the metric $D_\delta$ 
which do  not contain the origin, i.e.\ balls with radius $r\leq  D_\delta((\xi_1,0)\,;(0,0))$.
Our aim is to prove the following.

\begin{prop}\label{ppi3.4}
There are $\lambda_2>0$ and $\kappa_2\in\langle0,1]$  such that 
\begin{equation}
\int_{B_\Delta((\xi_1,0);r)}dx\,\Gamma_\delta(\varphi)(x)\geq\lambda_2\,r^{-2}\int_{B_\Delta((\xi_1,0);\kappa_2 r)}dx\,(\varphi(x)-\langle\varphi\rangle)^2
\label{epi3.34}
\end{equation}
for all $\varphi\in C_c^1(\Ri^{n+m})$  and $r\in\langle0,r_{\xi}\,]$ where $r_{\xi}=D_\delta((0,0)\,;(\xi_1,0))$ and $\langle\varphi\rangle$ is the average of $\varphi$ over $B_\Delta((\xi_1,0);\kappa_2 \,r)$.
\end{prop}

The proof has several features in common with Case~I. 
It relies in part on estimating on special sets which are are embedded in an appropriate manner in the Riemannian balls.
These sets are defined 
for  each $\xi_1\in\Ri^n$ and $\kappa\in \langle0,1]$ by
\[
C(\xi\,;\kappa)=\{(x_1,x_2):|x_1-\xi_1|< (\kappa/2)\,r_\xi^{(\alpha,\alpha')}\,,\,|x_2|< (\kappa/2)\,r_\xi^{(\beta,\beta')}\,\}
\;.
\]
Thus $C(\xi\,;\kappa)$ is the product of an $n$-dimensional Euclidean ball centred at $\xi_1$ and  an $m$-dimensional 
Euclidean ball centred at $0$, both with diameter $\kappa$, rescaled
 by the Riemannian shape factors $r_\xi^{(\alpha,\alpha')}$ and $r_\xi^{(\beta,\beta')}$.
 The choice of balls instead of cubes is for convenience in the following estimates and is of no great significance.

The Riemannian rescaling ensures the following embedding in analogy with Lemma~\ref{lpi3.3}.

\begin{lemma}\label{lpi3.5}
$\;\;\; C(\xi\,;\kappa)\subseteq B_\Delta((\xi_1,0)\,;\kappa\,r_\xi)$ for all $\kappa\in \langle0,1]$.
\end{lemma}
\proof\
First note that  $r_\xi=|\xi_1|^{(1-\delta_1,1-\delta_1')}$.
Secondly, if  $x\in C(\xi\,;\kappa)$ then
\begin{eqnarray*}
D_\delta((\xi_1,0)\,;(x_1,x_2))&<& |x_1-\xi_1|\,|\xi_1|^{(-\delta_1,-\delta_1')}
+|x_2|\,|\xi_1|^{(-\delta_2,-\delta_2')}
\\[5pt]
&< &(\kappa/2)\,r_\xi^{(\alpha,\alpha')}\,r_\xi^{(-\alpha\delta_1,-\alpha'\delta_1')}+
(\kappa/2)\,r_\xi^{(\beta,\beta')}\,r_\xi^{(-\delta_2\alpha,-\delta_2'\alpha')}\,
\\[5pt]
&=&(\kappa/2)\,r_\xi+(\kappa/2)\,r_\xi=\kappa\,r_\xi
\;.
\end{eqnarray*}
Thus $x\in B_\Delta((\xi_1,0)\,;\kappa\,r_\xi)$ and the embedding is established.
\hfill$\Box$
\bigskip

The starting point for the derivation of the Poincar\'e inequality for the balls $B_\Delta((\xi,0);r)$ is the following analogue of 
Proposition~\ref{ppi3.1}.

\begin{prop}\label{ppi3.5}
There is a $\lambda>0$ such that 
\begin{equation}
\int_{C(\xi;\kappa)}dx\,\Gamma_\delta(\varphi)(x)\geq\lambda\,(\kappa\,r_\xi)^{-2}\int_{C(\xi;\kappa)}dx\,(\varphi(x)-\langle\varphi\rangle)^2
\label{epi3.10}
\end{equation}
for all $\kappa\in\langle0,1]$ and all $\varphi\in C_c^1(\Ri^{n+m})$ where $\langle\varphi\rangle$ is the average of $\varphi$ over ${C(\xi;\kappa)}$.
\end{prop}
\proof\
Let  $x\in C(\xi;\kappa)$. 
Since $|x_1|\geq |\xi_1|-|x_1-\xi_1|$ and $|\xi_1|=r_\xi^{(\alpha, \alpha')}$ it follows that 
$|x_1|\geq r_\xi^{(\alpha, \alpha')}-(\kappa/2)\,r_\xi^{(\alpha,\alpha')}\geq 2^{-1}r_\xi^{(\alpha,\alpha')}$.
Therefore 
\[
|x_1|^{(2\delta_i,2\delta_i')}\geq (2^{-1}r_\xi^{(\alpha,\alpha')})^{(2\delta_i,2\delta_i')}\geq (1/8)\,
r_\xi^{(2\alpha\delta_i,2\alpha'\delta_i')}
\]
where the second estimate uses Proposition~\ref{ppi2.1}.
Consequently,
\[
\Gamma_\delta(\varphi)(x)\geq (1/8)\,\Big(r_\xi^{(2\alpha\delta_1,2\alpha'\delta_1')}((\nabla_{x_1}\varphi)(x))^2
+r_\xi^{(2\alpha\delta_2,2\alpha'\delta_2')}((\nabla_{x_2}\varphi(x))^2\Big)
\]
for all $x\in C(\xi;\kappa)$.

Next changing integration  variables to
 $y_1=(x_1-\xi_1)/r_\xi^{(\alpha,\alpha')}$ and $y_2=x_2/r_\xi^{(\beta,\beta')}$, one calculates that
 \begin{eqnarray*}
\int_{C(\xi;\kappa)}dx\,\Gamma_\delta(\varphi)(x)&\geq&
J\int_{B_\kappa}dy_1\int_{C_\kappa}dy_2\,\Big(
r_\xi^{(-2\alpha, -2\alpha')}r_\xi^{(2\delta_1\alpha, 2\delta_1'\alpha')}((\nabla_{y_1}\varphi)(y))^2\\[5pt]
&&\hspace{6cm}
+r_\xi^{(-2\beta, -2\beta')}r_\xi^{(2\delta_2\alpha, 2\delta_2'\alpha')}((\nabla_{y_2}\varphi)(y))^2\Big)\\[5pt]
&=&
J\,r_\xi^{-2}\int_{B_\kappa}dy_1\int_{C_\kappa}dy_2\,\Big(
((\nabla_{y_1}\varphi)(y))^2+((\nabla_{y_2}\varphi)(y))^2\Big)\\[5pt]
&=&J\,r_\xi^{-2}\int_{B_\kappa}dy_1\int_{C_\kappa}dy_2\,((\nabla_y\varphi)(y))^2
\end{eqnarray*}
 where $B_{\kappa}=\{y_1\in\Ri^n:|y_1|< (\kappa/2)\}$, $C_{\kappa}=\{y_2\in\Ri^m:|y_2|< (\kappa/2)\}$ and  $J$ is the Jacobian of the coordinate transformation.

Now one can use the usual Poincar\'e inequality for the Laplacian on the set $B_\kappa\times C_\kappa$ to deduce that there is a $\lambda>0$, independent of $\kappa$, such that
\begin{eqnarray*}
\int_{B_\kappa}dy_1\int_{C_\kappa}dy_2\,((\nabla_y\varphi)(y))^2
&\geq &\lambda\,\kappa^{-2}\int_{B_\kappa}dy_1\int_{C_\kappa}dy_2\,(\varphi(y)-\langle\varphi\rangle)^2
\\[5pt]
&=&\lambda\,\kappa^{-2}\inf_{M\in\Ri}\int_{B_\kappa}dy_1\int_{C_\kappa}dy_2\,(\varphi(y)-M)^2
\;.
\end{eqnarray*}
In particular  $\lambda$ is the lowest eigenvalue of the Laplacian on the set of
 $y\in\Ri^{n+m}$ with $|y_1|\leq 1/2$ and $|y_2|\leq 1/2$.
Consequently it is independent of all 
the parameters $\xi$, $\kappa$, $\delta_i$ etc.
Combining these estimates and reverting to the original coordinates one deduces that 
 \begin{eqnarray*}
\int_{C(\xi;\kappa)}dx\,\Gamma_\delta(\varphi)(x)&\geq &
J\,\lambda\,(\kappa\,r_\xi)^{-2}\inf_{M\in\Ri}\int_{B_\kappa}dy_1\int_{C_\kappa}dy_2\,(\varphi(y)-M)^2
\\[5pt]
&=&\lambda\,(\kappa\,r_\xi)^{-2} \inf_{M\in\Ri}\int_{C(\xi;\kappa)}dx\,(\varphi(x)-M)^2\\[5pt]
&=&\lambda\,(\kappa\,r_\xi)^{-2} \int_{C(\xi;\kappa)}dx\,(\varphi(x)-\langle\varphi\rangle)^2
\end{eqnarray*}
for all $\kappa\in\langle0,1]$ and all $\varphi\in C_c^1(\Ri^{n+m})$\hfill$\Box$

\bigskip

Next one can transfer the Poincar\'e inequality on the $C(\xi\,;\kappa)$ to the balls $B_\Delta((\xi_1,0)\,;r)$ by the  following  embedding analogous to Lemma~\ref{lpi3.4}.

\begin{lemma}\label{lpi3.6} 
There is a $\kappa_0\in\langle0, 1]$ such that $B_\Delta((\xi_1,0)\,;\kappa_0\kappa \,r_\xi)\subseteq C(\xi\,;\kappa)$ for all $\kappa\in \langle0,1]$.
\end{lemma}
\proof\
Consider the family of balls $B_\Delta((\xi_1,0)\,;\kappa\,r_\xi)$ for $\kappa\in\langle0,1]$
and introduce the set  $B_n=\{x_1\in \Ri^n: (x_1,0)\in B_\Delta((\xi_1,0)\,;\kappa\,r_\xi)\}$.
Then $x_1\in B_n$ if and only if
\[
|x_1-\xi_1|< \kappa\,r_\xi\,(|x_1|+|\xi_1|)^{(\delta_1,\delta_1')}
\]
i.e.\  $D_\delta((x_1,0)\,;(\xi_1,0))<\kappa\,r_\xi$.
Therefore $B_n\subseteq C_n$ where
\begin{eqnarray*}
C_n&=&\{x_1\in \Ri^n: |x_1-\xi_1|< \kappa\,r_\xi\,(|x_1-\xi_1|+2\,|\xi_1|)^{(\delta_1,\delta_1')}\}\\[5pt]
&=&\{x_1\in \Ri^n: |x_1-\xi_1|< \kappa\,|\xi_1|^{(1-\delta_1,1-\delta_1')}\,(|x_1-\xi_1|+2\,|\xi_1|)^{(\delta_1,\delta_1')}\}
\;.
\end{eqnarray*}
Since $|\xi_1|\neq0$ one deduces that $x_1\in C_n$ if and only if 
\[
|x_1-\xi_1|/|\xi_1|<\kappa\,|\xi_1|^{(-\delta_1,-\delta_1')}\,\Big(|\xi_1|(2+|x_1-\xi_1|/|\xi_1|)\Big)^{(\delta_1,\delta_1')}
\;.
\]
Therefore by Proposition~\ref{ppi2.1}  it is necessary that 
\[
|x_1-\xi_1|/|\xi_1|< 4\,\kappa\,(2+|x_1-\xi_1|/|\xi_1|)^{\delta_1\vee\delta_1'}
\;.
\]
Then setting $\sigma=|x_1-\xi_1|/(\kappa\,|\xi_1|)$ one must have
$\sigma\leq 4\,(2+\sigma)^{\delta_1\vee\delta_1'}$.
But $\delta_1\vee\delta_1'<1$ and one concludes that $\sigma\leq a_1$ where  $a_1>0$ is the unique solution of 
$a_1=4\,(2+a_1)^{\delta_1\vee\delta_1'}$.
Consequently
\[
|x_1-\xi_1|\leq a_1\,\kappa\,|\xi_1|=a_1\,\kappa\,r_\xi^{(\alpha, \alpha')}
\]
for all $\kappa\in\langle0,1]$.

Next observe that  if  $x\in B_\Delta((\xi_1,0)\,;\kappa\,r_\xi)$ then $x_1\in B_n$.
This follows by contradiction.
Assume  $x\in B_\Delta((\xi_1,0)\,;\kappa\,r_\xi)$ but $x_1\not\in B_n$.
Then
\begin{eqnarray*}
\kappa\,r_\xi> D_\delta((x_1,x_2)\,;(\xi_1,0))\geq |x_1-\xi_1|/(|x_1|+|\xi|_1)^{(\delta_1,\delta_1')}
=D_\delta((x_1,0)\,;(\xi_1,0))\geq\kappa\,r_\xi
\;.
\end{eqnarray*}
The first inequality follows since $x\in B_\Delta((\xi_1,0)\,;\kappa\,r_\xi)$  and the last follows since $x_1\not\in B_n$.
But  this gives a contradiction.

It now  follows that  there is a  $\rho>0$, dependent on $\xi$, such that 
\[
\inf_{x_1\in B_n, |x_2|\geq\rho}D_\delta((x_1,x_2)\,;(\xi_1,0))=\kappa\,r_\xi
\;.
\]
Then $B_\Delta((\xi_1,0)\,;\kappa\,r_\xi)\subseteq B_n\times\{x_2\in \Ri^m: |x_2|<\rho\}$.
Next we estimate $\rho$.

First it follows from the observation
\[
\kappa\,r_\xi=\inf_{x_1\in B_n, |x_2|\geq\rho} D_\delta((x_1,x_2)\,;(\xi_1,0))
\leq \inf_{ |x_2|\geq\rho}D_\delta((\xi_1,x_2)\,;(\xi_1,0))
\]
that
\[
\kappa\,r_\xi \leq \rho\left(|\xi_1|^{(\delta_2,\delta_2')}+\rho^{(\gamma,\gamma')}\right)^{-1}
\;.
\]
Therefore one obtains a lower bound on $\rho$,
\begin{eqnarray*}
\rho\geq \kappa\,r_\xi\,|\xi_1|^{(\delta_2,\delta_2')}\geq \kappa\,r_\xi^{(\beta,\beta')}
\end{eqnarray*}
where the  second step uses $r_\xi=|\xi_1|^{(1-\delta_1,1-\delta_1')}$.

Secondly, one obtains an upper  bound on $\rho$ by observing that 
\begin{eqnarray*}
\kappa\,r_\xi &\geq&
\inf_{x_1\in B_n, |x_2|\geq\rho}|x_2|\left((|x_1|+|\xi_1|)^{(\delta_2,\delta_2')}+(|x_2|)^{(\gamma,\gamma')}\right)^{-1}\\[5pt]
&\geq &\inf_{x_1\in B_n}\rho\left((|x_1|+|\xi_1|)^{(\delta_2,\delta_2')}+\rho^{(\gamma,\gamma')}\right)^{-1}
\;.
\end{eqnarray*}
Since by the previous estimate $x_1\in B_n$  satisfies the bound $|x_1|\leq (1+a\,\kappa)\,|\xi_1|\leq (1+a)\,|\xi_1|$ it follows that 
\begin{eqnarray*}
\rho&\leq& \kappa\,r_\xi\,\left(((1+a)\,|\xi_1|)^{(\delta_2,\delta_2')}+\rho^{(\gamma,\gamma')}\right)
\end{eqnarray*}
for all $\kappa\in\langle0,1]$.
Then by   Proposition~\ref{ppi2.1} one deduces that there is a $b>0$ such that
\[
\rho\leq \kappa\,r_\xi\,\left(b\,|\xi_1|^{(\delta_2,\delta_2')}+\rho^{(\gamma,\gamma')}\right)
\]
for all $\kappa\in\langle0,1]$.
But $r_\xi= |\xi_1|^{(1-\delta_1,1-\delta_1')}$ so
\begin{equation}
\rho\leq 
b\,\kappa\,|\xi_1|^{(\tau,\tau')}+\kappa\,|\xi_1|^{(1-\delta_1,1-\delta_1')}\rho^{(\gamma,\gamma')}
\;.
\label{e4.5}
\end{equation}
with  $\tau=1+\delta_2-\delta_1$ and $\tau'=1+\delta_2'-\delta_1'$.
Now we estimate $\rho$ in two separate cases, $|\xi_1|\leq1$ and $|\xi_1|\geq1$.

If $|\xi_1|\leq 1$ 
then 
\[
\rho\leq 
b\,\kappa+\kappa\,\rho^{(\gamma,\gamma')}
\;.
\]
Since $\gamma, \gamma'<1$ it follows that there is a $b_0>0$ such that $\rho\leq b_0$.
Now if $b_0\leq 1$ then $\rho\leq (1+b)\,\kappa$.
Alternatively if $b_0\geq1$ then $b_0/\kappa\leq b+b_0^{\gamma'}\leq b+(b_0/\kappa)^{\gamma'}$ and 
$b_0/\kappa$ must be uniformly bounded.
Therefore there is a $b_0>0$ such that $\rho\leq b_0\,\kappa$ for all $\kappa\in\langle0,1]$.
Hence by another application of Proposition~\ref{ppi2.1} one deduces that there is a $c>0$ such that 
$\rho^{(\gamma,\gamma')}\leq c\,(\rho/(b\,\kappa))^{\gamma}$ for all $\kappa\in\langle0, 1]$.
The value of $c$ is independent of $\kappa$ and $\xi_1$.
Then it follows from (\ref{e4.5}) that 
\begin{eqnarray*}
\rho/(b_0\,\kappa\,|\xi_1|^\tau)&\leq& (b/b_0)+(c/b_0)\,|\xi_1|^{1-\delta_1-\tau}\,(\rho/(b_0\,\kappa))^{\gamma}\\[5pt]
&=&(b/b_0)+(c/b_0)\,(\rho/(b_0\,\kappa\,|\xi_1|^\tau))^{\gamma}
\;.
\end{eqnarray*}
Therefore, since $\gamma<1$ one has 
\[
\rho\leq a_2\,\kappa\,|\xi_1|^\tau=a_2\,\kappa\,r_\xi^\beta
\]
where $a_2$ satisfies $a_2=b+c\,(a_2/b_0)^{\gamma}$.

Next suppose $|\xi_1|\geq 1$.
It  then follows from the above discussion of the lower bound on $\rho$ that 
\[
 \rho\geq \kappa\,r_\xi\,|\xi_1|^{\delta_2'}
 \geq \kappa\,|\xi_1|^{\tau'}
 =\kappa\,r_\xi^{\beta'}
 \]
 because $r_\xi=|\xi_1|^{1-\delta_1'}$.
 In particular $\rho/\kappa\geq1$.
 But then
 \[
 \rho^{(\gamma,\gamma')}=(\kappa\,(\rho/\kappa))^{(\gamma,\gamma')}
 \leq 2^{2(\gamma\vee\gamma')}\,\kappa^{(\gamma\wedge\gamma')}\,(\rho/\kappa)^{\gamma'}
 \]
 by Proposition~\ref{ppi2.1}.
 Hence it follows from (\ref{e4.5}) that 
\begin{eqnarray*}
 (\rho/(\kappa\,|\xi_1|^{\tau'}))&\leq& 
b+4\,|\xi_1|^{-\tau'}\,|\xi_1|^{1-\delta_1'}\,\kappa^{(\gamma\wedge\gamma')}\,(\rho/\kappa)^{\gamma'}\\[5pt]
&\leq&b+4\, (\rho/(\kappa\,|\xi_1|^{\tau'}))^{\gamma'}
\end{eqnarray*}
where  the second step uses $\kappa\leq 1$.
Consequently one deduces as before that 
\[
\rho \leq a_2\,\kappa\,|\xi_1|^{\tau'}= a_2\,\kappa\,r_\xi^{\beta'}
\]
for all $\kappa\in\langle0,1]$ and $|\xi_1|\geq1$
where $a_2=b+4\,a_2^{\gamma'}$.

Finally combination of these results leads to the conclusion that there are
$a_2,a_2'>0$ such that 
\[
a_2'\,\kappa\,|\xi_1|^{(\tau,\tau')}\leq \rho\leq a_2 \,\kappa\,|\xi_1|^{(\tau,\tau')}
\]
or, equivalently,
\[
a_2'\,\kappa\,r_\xi^{(\beta,\beta')}\leq \rho\leq a_2\,\kappa\,r_\xi^{(\beta,\beta')}
\]
for all $\kappa\in\langle0,1]$.
The values of $a_2$ and $a_2'$ are independent of $\xi_1$ and $\kappa$.

\smallskip

Now we can complete the proof of Lemma~\ref{lpi3.6}.

\smallskip

If $x=(x_1,x_2)\in B_\Delta((\xi_1,0)\,;\kappa_0\,\kappa\,r_\xi)$ the foregoing estimates are valid with
$\kappa$ replaced by $\kappa_0\,\kappa$ with $\kappa_0\in\langle0,1]$ and   $\kappa\in\langle0,1]$.
Therefore
\[
|x_1-\xi_1|\leq  a_1\,\kappa_0\,\kappa\,r_\xi^{(\alpha, \alpha')}\;\;\;\;{\rm and}\;\;\;\;\;
|x_2|\leq a_2\,\kappa_0\,\kappa\,r^{(\beta,\beta')}
\;.
\]
Hence if $(a_1\vee a_2)\,\kappa_0< 1/2$ it follows that $x\in C(\xi\,;\kappa)$.
\hfill$\Box$

\bigskip

At this point the proof of Proposition~\ref{ppi3.4} is immediate.
First there is a $\lambda>0$ such that 
\begin{eqnarray*}
\int_{B_\Delta((\xi_1,0);\kappa r_\xi)}dx\,\Gamma_\delta(\varphi)(x)&\geq& \int_{C(\xi;\kappa)}dx\,\Gamma_\delta(\varphi)(x)\\[5pt]
&\geq&\lambda\,(\kappa\,r_\xi)^{-2}\int_{C(\xi;\kappa)}dx\,(\varphi(x)-\langle\varphi\rangle)^2
\end{eqnarray*}
for all $\kappa\in\langle0,1]$ by Lemma~\ref{lpi3.5} and Proposition~\ref{ppi3.4}.

Secondly, there is a $\kappa_0\in\langle0,1]$ such that 
\begin{eqnarray*}
\int_{C(\xi;\kappa)}dx\,(\varphi(x)-\langle\varphi\rangle)^2&=&\inf_{M\in\Ri}\int_{C(\xi;\kappa)}dx\,(\varphi(x)-M)^2\\[5pt]
&\geq&\inf_{M\in\Ri}\int_{B_\Delta((\xi_1,0);\kappa_0\kappa\,r_\xi)}dx\,(\varphi(x)-M)^2\\[5pt]
&=&\int_{B_\Delta((\xi_1,0);\kappa_0\kappa\,r_\xi)}dx\,(\varphi(x)-\langle\varphi\rangle)^2
\end{eqnarray*}
for all $\kappa\in\langle0,1]$ by two more applications of (\ref{epi2.0}) and by Lemma~\ref{lpi3.6}.

Therefore one concludes that 
 \begin{eqnarray*}
 \int_{B_\Delta((\xi_1,0);\kappa\,r_\xi)}dx\,\Gamma_\delta(\varphi)(x)
 &\geq&\lambda\,(\kappa\,r_\xi)^{-2} \int_{B_\Delta((\xi_1,0);\kappa_0\kappa\,r_\xi)}dx\,(\varphi(x)-\langle\varphi\rangle)^2
 \end{eqnarray*}
 for all $\kappa\in\langle0,1]$ and the statement of Proposition~\ref{ppi3.4}, with $\lambda_2=\lambda$ and $\kappa_2=\kappa_0$,  follows by setting $r=\kappa\,r_\xi$.
 \hfill$\Box$

\bigskip

Thus Proposition~\ref{ppi3.4}  establishes the Poincar\'e inequality for the  balls $B_\Delta((\xi_1,0)\,;r)$ for $\xi\neq0$ and for  all $r\leq D_\delta((\xi_1,0)\,;(0,0))$.

\begin{remarkn}\label{rpi3.1}
Note that in the foregoing proof of Case~II we do not need to assume that $\delta_1\vee\delta_1'\in[0,1/2\rangle$  if $n=1$.
\end{remarkn}
\bigskip

\noindent{\bf Case III--General balls.}
To complete the proof of the Poincar\'e inequality (\ref{epi1.5}) it suffices to verify it for the balls $B_\Delta((\xi_1,0)\,;r)$ with $\xi_1\neq0$ and  $r\geq r_\xi$ where again $r_\xi=D_\delta((0,0)\,;(\xi_1,0))$.
(If $\xi_1=0$ the inequality follows for all $r>0$ by Case~I and if  $r\leq r_\xi$ then it follows from Case~II.)
The general case is  a corollary of the two special cases. 

First assume $r\geq K\,r_\xi$ with $K=2\,(1+\kappa_1)/\kappa_1$ where $\kappa_1$ is the parameter of Proposition~\ref{ppi3.3}.
Then $r\geq 2\,r_\xi$ and 
 $B_\Delta((0,0)\,;r-r_\xi)\subseteq B_\Delta((\xi_1,0)\,; r)$.
 Therefore 
 \[
 \int_{B_\Delta((\xi_1,0); r)}\Gamma(\varphi)\geq \int_{B_\Delta((0,0);r-r_\xi)}\Gamma(\varphi)
 \geq \lambda\,(r-r_\xi)^{-2}\inf_{M\in \Ri}\int_{B_\Delta((0,0);\kappa_1(r-r_\xi))}(\varphi-M)^2
 \]
 by Proposition~\ref{ppi3.3}.
 But $0<r-r_\xi<r$. Hence $(r-r_\xi)^{-2}>r^{-2}$.
 Moreover, one has the inclusion $B((\xi_1,0)\,;\kappa_1(r-r_\xi)-r_\xi))\subseteq B_\delta((0,0)\,;\kappa_1(r-r_\xi))$.
 Therefore 
 \[
 \int_{B_\Delta((\xi_1,0); r)}\Gamma(\varphi)\geq  \lambda\,r^{-2}\inf_{M\in \Ri}\int_{B_\Delta((\xi_1,0);\kappa_1(r-r_\xi)-r_\xi)}(\varphi-M)^2
 \;.
 \]
Since $\kappa_1(r-r_\xi)-r_\xi\geq \kappa_1\,r/2$ it then follows that 
\[
 \int_{B_\Delta((\xi_1,0); r)}\Gamma(\varphi)\geq  \lambda\,r^{-2}\inf_{M\in \Ri}\int_{B_\Delta((\xi_1,0);\kappa_1 r/2)}(\varphi-M)^2
 \]
for all $r\geq K\,r_\xi$.

Secondly suppose $K\,r_\xi\geq r\geq r_\xi$.
Then
\[
\int_{B_\Delta((\xi_1,0); r)}\Gamma(\varphi)\geq 
\int_{B_\Delta((\xi_1,0); r_\xi)}\Gamma(\varphi)\geq \lambda\,(r_\xi)^{-2}
\inf_{M\in \Ri}\int_{B_\Delta((\xi_1,0);\kappa_2 r_\xi)}(\varphi-M)^2
\]
by Proposition~\ref{ppi3.4}.
But $\kappa_2r_\xi\geq (\kappa_2/K)r$ and $(r_\xi)^{-2}\geq r^{-2}$.
Therefore
\[
\int_{B_\Delta((\xi_1,0); r)}\Gamma(\varphi)\geq \lambda\, r^{-2}\inf_{M\in \Ri}\int_{B_\Delta((\xi_1,0);(\kappa_2/K)r)}(\varphi-M)^2
\]
for all $r\in [r_\xi, 2(1+\kappa_1)\,r_\xi/\kappa_1]$.

Combination of these two results  establishes the Poincar\'e inequality for the balls $B_\Delta((\xi_1,0)\,;r)$ and  all $r\geq r_\xi$ with the value of $\kappa$ given by $(\kappa_2/K)\wedge (\kappa_1/2)$.
But  the bounds were established for $\xi_1=0$ and all $r$ in Case~I and for $\xi_1\neq0$ and  $r\leq r_\xi$ in Case~II.
Thus  it follows that the Poincar\'e inequality is valid
 for the $B_\Delta((\xi_1,0)\,;r)$ for all $\xi_1\in \Ri^n$ and all $r>0$.
 Finally, since the Riemannian metric $d(\,\cdot\,;\,\cdot\,)$ is equivalent to the metric $D_\delta(\,\cdot\,;\,\cdot\,)$ 
which defines the balls $B_\Delta(0\,;r)$ the Poincar\'e inequality is valid for the Riemannian balls
 by the discussion in Section~\ref{S2}.
The change to  an equivalent metric only requires a change in the value of the parameter $\kappa$ in the inequality.

\medskip

At this stage we have established the first statement of Theorem~\ref{tpi1.1}.
Next we prove the second statement, the failure of the Poincar\'e inequality for $n=1$ and  $\delta_1\vee\delta_1'\in[1/2,1\rangle$.
We will establish this in two steps.
\medskip

First assume $\delta_1\in [1/2,1\rangle$. 
Let $\varphi=\chi\psi$ where $\chi\in C_c^1(\Ri)$, $\psi\in C_c^1(\Ri^m)$ and $\psi$ is equal to one on $[-1,1]^m$.
Then
\[
\int_{[-1,1]^{1+m}}dx\,\Gamma_\delta(\varphi)(x)=2^m\int^1_{-1}dx_1 \,\Gamma_\delta^{(1)}(\chi)(x_1)
\]
where $\Gamma_\delta^{(1)}$ is the  {\it carr\'e du champ} of 
$H_\delta^{(1)}=-d_{x_1}\,|x_1|^{2\delta_1}\,d_{x_1}$ acting on $L_2(-1,1)$. 
Explicitly
$\Gamma_\delta^{(1)}(\chi(x_1))=|x_1|^{2\delta_1}(\chi'(x_1))^2$.
Moreover,
\[
\int_{[-1,1]^{1+m}}dx\,(\varphi(x)-\langle\varphi\rangle)^2=
2^m\int^1_{-1}dx_1\,(\chi(x_1)-\langle\chi\rangle)^2
\]
where $\langle\varphi\rangle$ is the average of $\varphi$ over  the cube $[-1,1]^{1+m}$ and $\langle\chi\rangle$ is the average of $\chi$ over the interval $[-1,1]$.
Therefore the Poincar\'e inequality fails for $H_\delta$ on $[-1,1]^{1+m}$ if it fails for $H_\delta^{(1)}$ on $[-1,1]$.

Define $\chi_n \colon \Ri \to [-1,1]$ by
\[
\chi_n(x)
= \left\{ \begin{array}{ll}
   0 & \mbox{if } 0\leq x\leq n^{-1} \\[5pt]
  1- \eta_n^{-1} \, \eta(x) & \mbox{if } n^{-1}\leq x\leq 1 \\[5pt]
   1 & \mbox{if } x \geq 1
          \end{array} \right.
          \]
  and
  \[
  \hspace{1.2cm}
\chi_n(x)
= \left\{ \begin{array}{ll}
   0 & \mbox{if }  -n^{-1} \leq x\leq 0  \\[5pt]
  -1+\sigma_n^{-1} \, \sigma(x) & \mbox{if }  -1\leq x\leq -n^{-1} \\[5pt]
   -1 & \mbox{if } x \leq -1 
          \end{array} \right.
         \]
where
\[
\eta(x)=\int^1_xds\, |s|^{-2\delta_1},\;\;
\sigma(x)=\int^x_{-1} ds\, |s|^{-2\delta_1},\;\;\eta_n=\eta(n^{-1}) \;\;{\rm and}\;
\sigma_n=\sigma(-n^{-1}) 
\;.
\]
Note that $\chi_n$ is an absolutely continuous increasing function.
Moreover, since $\delta_1\in[1/2,1\rangle$, it follows that $\lim_{n \to \infty} \chi_n(x) = 1$ if $x>0$ and  $\lim_{n \to \infty} \chi_n(x) = -1$ if $x<0$.
For example if $\delta_1=1/2$ then $\eta(x)\sim \log |x|\sim \sigma(x)$
and  $\eta_n\sim \log n\sim\sigma_n$.
Therefore $\langle\chi_n\rangle\to0$ as $n\to\infty$ and
\[
\lim_{n\to\infty}\int^1_{-1} dx\,\left(\chi_n(x)-\langle\chi_n\rangle\right)^2=
\lim_{n\to\infty}\int^1_{-1} dx\,\chi_n(x)^2=2
\;.
\]
But 
\[
\lim_{n\to\infty}\int^1_{-1}dx\, |x|^{2\delta_1}\,(\chi_n'(x))^2=\lim_{n\to\infty}(\sigma_n+\eta_n)=0
\;.
\]
Therefore the Poincar\'e inequality for $H_\delta^{(1)}$ on $[-1,1]$ must fail for $\chi_n$ if $n$ is sufficiently large.
Consequently the  Poincar\'e inequality (\ref{epi3.1}) for $H_\delta$  must fail for $\varphi=\chi_n\psi$ for large $n$.

\medskip

Secondly, assume $\delta_1\in [0,1/2\rangle$ and  $\delta_1'\in [1/2,1\rangle$. 
We aim to show that the Poincar\'e inequality fails for Riemannian balls of large radius centred at the origin.
But it follows from the discussion of Case~I above that it suffices to prove that (\ref{epi3.22}) fails for centred cubes
$C_t$ with $t$ large.
This can again be reduced to a one-dimensional problem.

For each $t>0$ set $\varphi=\chi\,\psi$ with $\chi\in C_c^1(\Ri)$ and $\psi\in C_c^1(\Ri^m)$ with $\psi(x_2)=1$ if 
$|x_2|_\infty<t^{(\beta,\beta')}$.
Then one has
\[
\int_{C_t}dx\,\Gamma_\delta(\varphi)(x)=2^m\int_{|x_1|<t^{(\alpha,\alpha')}}dx_1 \,\Gamma_\delta^{(1)}(\chi)(x_1)
\]
and 
\[
\int_{C_t}dx\,(\varphi(x)-\langle\varphi\rangle)^2=
2^m\int_{|x_1|<t^{(\alpha,\alpha')}}dx_1\,(\chi(x_1)-\langle\chi\rangle)^2
\]
where $\langle\varphi\rangle$ is the average of $\varphi$ over  the cube $C_t$ and $\langle\chi\rangle$ is the average of $\chi$ over the interval $I_t=\{x_1\in\Ri:|x_1|<t^{(\alpha,\alpha')}\}$.
Thus to establish that (\ref{epi3.22}) fails for $C_t$ it suffices to establish that the one-dimensional analogue fails on $I_t$.

First we consider the particular case $\delta_1\in [0,1/2\rangle$ but  $\delta_1'=1/2$.
Then $\alpha\in [1,2\rangle$ but  $\alpha'=2$.
Let $\chi$ be an odd function with $\chi(x_1)=\int^{x_1}_0ds\,s^{(-2\delta_1,-1)}$ for $x_1\geq0$.
Then $\chi$ is locally bounded  and $\chi(x_1)\sim\pm\log|x_1|$ as $x_1\to\pm\infty$.
Now $\Gamma_\delta^{(1)}(\chi)(x_1)=|x_1|^{(-2\delta_1,-1)}$ and it follows that 
\[
\int_{|x_1|<t^{(\alpha,2)}}dx_1 \,\Gamma_\delta^{(1)}(\chi)(x_1)\sim \int^{t^2}_1ds\,s^{-1}\sim \log t
\]
as $t\to \infty$.
But  $\langle\chi\rangle=0$ because the function is odd and 
\[
t^{-2}\int_{|x_1|<t^{(\alpha,2)}}dx_1\,(\chi(x_1)-\langle\chi\rangle)^2
=t^{-2}\int_{|x_1|<t^{(\alpha,2)}}dx_1\,(\chi(x_1))^2
\sim\int^{t^2}_1ds\,(\log s)^2
\sim (\log t)^2
\]
as $t\to\infty$.
Thus the Poincar\'e inequality must fail for large $t$.

Secondly consider the  case $\delta_1\in [0,1/2\rangle$ but  $\delta_1'\in \langle1/2,1\rangle$.
Again $\alpha\in [1,2\rangle$ but now $\alpha'>2$.
Let $\chi\in C^1(\Ri)$ be an odd  increasing function  with $\chi(x_1)=1$ if $x_1\geq 1$.
Then $\Gamma_\delta^{(1)}(\chi)$ is a positive  bounded function with support in the interval $[-1,1]$.
Hence $\int_{|x_1|<t^{(\alpha,\alpha')}}dx_1 \,\Gamma_\delta^{(1)}(\chi)(x_1)$ is bounded uniformly for $t\geq1$.
On the other hand 
\[
t^{-2}\int_{|x_1|<t^{(\alpha,\alpha')}}dx_1\,(\chi(x_1)-\langle\chi\rangle)^2
=t^{-2}\int_{|x_1|<t^{(\alpha,\alpha')}}dx_1\,(\chi(x_1))^2\sim t^{-2}t^{\alpha'}
\]
as $t\to\infty$.
Since $\alpha'>2$ the Poincar\'e inequality must again fail.
\medskip

These examples establish the second statement of Theorem~\ref{tpi1.1} 
and it remains to prove the third statement.

\bigskip

The proof is by modification of the above argument for the Poincar\'e inequality on $\Ri^{n+m}$.
The only significant modification occurs in the discussion of the (half) balls centred at the origin.
Consider the case of $\Ri_+\times \Ri^m$.
Then $B_+(0\,;r)=B(0\,;r)\cap \{x_1:x_1>0\}$.
The proof of (\ref{epi1.6}) for $B_+(0\,;r)$ begins with the analogue of Proposition~\ref{ppi3.1}.

\begin{prop}\label{ppi3.6}
Assume $\delta_1\in[1/2,1\rangle$.
Then there is a $\lambda>0$ such that 
\begin{equation}
\int^1_0dx_1\int_{[-1,1]^{m}}dx_2\,\Gamma_\delta(\varphi)(x_1,x_2)\geq\lambda
\int^1_0dx_1\int_{[-1,1]^{m}}dx_2\,(\varphi(x_1,x_2)-\langle\varphi\rangle)^2
\label{epi3.11}
\end{equation}
for all $\varphi\in C_c^1(\Ri_+\times\Ri^{m})$.
\end{prop}
\proof\ The argument used to prove Proposition~\ref{ppi3.1} is easily adapted to the half-space and again
reduces the problem to a pair of one-dimensional problems.
It is reduced to proving (\ref{epi3.2}) and (\ref{epi3.3}) with the interval $[-1,1]$ replaced by $[0,1]$ where $\tilde h_\delta$
is given by (\ref{epi3.4}) modified similarly.
But (\ref{epi3.2})  follows because $\psi\in C_c^1(\Ri)\mapsto \int^1_0dx\,x^{2\delta_1}(\psi'(x))^2$ is a closable form and its closure corresponds to the self-adjoint extension of the operator $-d_x\,x^{2\delta_1}\,d_x$ on $L_2(0,1)$ with Neumann boundary conditions at each endpoint.
This operator has, however, a compact resolvent and the lowest eigenvalue is zero with the constant function one as corresponding eigenvalue.
Since the condition $\int^1_0dx\,x^{2\delta_1}(\psi'(x))^2=0$ implies that $\psi'=0$ and  $\psi$ is constant
it follows that the zero eigenvalue is simple.
Thus (\ref{epi3.2})  is satisfied with $\lambda$  the second eigenvalue.
Note that this argument, in contrast to that used to prove Lemma~\ref{lpi3.1},  does not require $\delta_1<1/2$ but works equally well for all $\delta_1\in[0,1\rangle$.
The point is that it is for   the operator on $[0,1]$.
Next the proof of (\ref{epi3.3}) as given in Lemma~\ref{lpi3.2} remains unchanged.
It is based on Lemma~\ref{lpi3.1} and hence it also does not require $\delta_1<1/2$ but is valid  for all $\delta_1\in [0,1\rangle$.\hfill$\Box$

\bigskip

The second step in the proof is an analogue of Proposition~\ref{ppi3.2}.

\begin{prop}\label{ppi3.7}
Assume $\delta_1\in[1/2,1\rangle$.
Then there is a $\lambda>0$ such that 
\[
\int_{C^\pm_t}dx\,\Gamma_\delta(\varphi)(x)\geq\lambda\,t^{-2}\int_{C^\pm_t}dx\,(\varphi(x)-\langle\varphi\rangle_\pm)^2
\]
for all $\varphi\in C_c^1(\Ri^{1+m})$  and $t>0$ where $C^+_t=C_t\cap\{x_1>0\}$, $C^-_t=C_t\cap\{x_1<0\}$
$\langle\varphi\rangle_\pm=|C^\pm_t|^{-1}\int_{C^\pm_t}dx\,\varphi(x)$.
\end{prop}
\proof\
The proof is identical to the proof of Proposition~\ref{ppi3.2} but is now based on Proposition~\ref{ppi3.6}
and scaling in the half-space.
The key point is that Proposition~\ref{ppi3.2} has to be applicable to the operator on the half space with coefficients
$|x_1|^{2(\delta_i\vee\delta_1', \delta_i\wedge\delta_i')}$.
This, however, only requires $\delta_1\vee\delta_1'\in[1/2,1\rangle$ and this is ensured if $\delta_1\in[1/2,1\rangle$.
There is no restraint on $\delta_1'$ it can take all values in $[0,1\rangle$.
\hfill$\Box$

\bigskip

The rest of the proof of the Poincar\'e inequality (\ref{epi1.6}) now follows by the argument used earlier 
to establish (\ref{epi1.5}).
The inequality for half balls centred at the origin follows from Proposition~\ref{ppi3.7} by slight modification of the earlier embedding arguments for cubes and balls. 
Then the proof for balls completely contained in the appropriate half-space follows by the discussion of Case~II in 
of the proof of  (\ref{epi1.5}).
This did not require the condition $\delta_1\vee\delta_1'\in[0,1/2\rangle$ (see Remark~\ref{rpi3.1})
 and applies equally well to the current situation with $\delta_1\in[1/2,1\rangle$ and $\delta_1'\in[0,1\rangle$.
Finally the inequality for general `balls' $B_\pm(\xi\,;r)$ follows as in the argument of Case~III above.
\hfill$\Box$

\section{Heat kernel bounds}\label{S4}

In this section we establish Theorems~\ref{tpi1.2} and \ref{tpi1.3}.
The upper Gaussian bounds follow from Corollary~6.6 of \cite{RSi2} once one establishes continuity of the kernel.
Therefore it suffices to prove the continuity and the lower Gaussian bounds.
These results are indirect  corollaries of Statements I and III of Theorem~\ref{tpi1.1}.
The key observation of Grigor'yan \cite{Gri4} and Saloff-Coste \cite{Sal2} is that the Poincar\'e inequality (\ref{epi1.5}) combined with the volume doubling property
of the Riemannian metric, \cite{RSi2} Corollary~5.2, implies the parabolic Harnack inequality of Moser \cite{Mos} 
on $\Ri^{n+m}$.
Similarly (\ref{epi1.6}) and volume doubling imply the Harnack inequality on $\Omega_\pm$.

The operator $H$ is defined to satisfy the (global) parabolic Harnack inequality on $\Ri^{n+m}$ if there exists an $a>0$
 such that for any $x\in\Ri^{n+m}$ and $t>0$ any non-negative (weak) solution $\varphi$ of the parabolic equation $(\partial_t+H)\varphi=0$ in the cylinder $Q=\langle t,\, t+r^2\rangle\times B(x\,;2\,r)$
satisfies
\begin{equation}
\sup_{Q^-}\varphi\leq a\,\inf_{Q^+}\varphi
\label{epi4.1}
\end{equation}
where $Q^-=[\,t+r^2\!/4,\, t+r^2\!/2\,]\times B(x\,;r)$ and $Q^+=[\, t+3\hspace{1pt} r^2\!/4, \,t+r^2\rangle\times B(x\,;r)$.
This definition is the key to establishing the continuity of the heat kernel and the lower Gaussian bounds.

\smallskip

\noindent{\bf Proof of Theorem~\ref{tpi1.2} }$\;$ 
First it follows from Theorem~\ref{tpi1.1}.I that the Poincar\'e inequality (\ref{epi1.5}) is valid.

Secondly, it follows from \cite{RSi2} Corollary~5.2 that the Riemannian balls
$B(x\,;r)$ satisfy the volume doubling property, i.e.\ there is a $b>0$ such that
\begin{equation}
|B(x\,;2\,r)|\leq b\,|B(x\,;r)|
\label{epi4.2}
\end{equation}
for all $x\in \Ri^{n+m}$  and all $r>0$.

Thirdly, Theorem~3.1 of \cite{Sal2} establishes that (\ref{epi1.5}) together with (\ref{epi4.2}) imply  that $H$ satisfies the parabolic Harnack inequality (\ref{epi4.2}).
Then, however, a straightforward argument of Moser,  \cite{Mos0} Section~5 or \cite{Mos} pages~108--109, establishes that each  non-negative solution $\varphi$ of $(\partial_t+H)\varphi=0$  is H\"older continuous.
Hence one deduces that the heat kernel $x,y\mapsto K_t(x\,;y)$ is jointly H\"older continuous.
Then the Gaussian upper bounds follow from Corollary~6.6 of \cite{RSi2}.
Moreover, the continuity ensures that the kernel $K_t$ is well-defined on the diagonal $x=y$
and it follows from Corollary~6.7 and Remark~6.8 of \cite{RSi2} that  there is a $c>0$ such that  \begin{equation}
K_t(x\,;x)\geq c\,|B(x\,;t^{1/2})|^{-1}
\label{epi4.4}
\end{equation}
for all $x\in \Ri^{n+m}$  and $t>0$.

Fourthly, fix $x$ and define
 define $\varphi$ by $\varphi(t,y)=K_t(x\,;y)$ for  all $t>0$.
Then  $\varphi$ is a non-negative weak solution of $(\partial_t+H)\varphi=0$ in the cylinder 
$Q=\langle 0,\,r^2\rangle\times B(x\,;2\,r)$.
Now if 
$Q^-=[\, r^2\!/4,r^2\!/2]\times B(x\,;r)$  and $Q^+= [3\hspace{1pt} r^2\!/4, \,r^2\,\rangle\times B(x\,;r)$ 
the parabolic Harnack inequality gives
\[
K_{r^2\!/2}(x\,;x)\leq \sup_{(t,y)\in Q^-}\varphi(t, y)
\leq a\,\inf_{(t,y)\in Q^+}\varphi(t,y)
\leq a \,K_{r^2}(x\,;y)\]
for all $y\in B(x\,;r)$.
Therefore, setting $r^2=t$, this estimate combined with  (\ref{epi4.4})  gives
\begin{equation}
K_t(x\,;y)\geq (c/a)\,|B(x\,;t^{1/2})|^{-1}
\label{epi4.41}
\end{equation}
for all $x\in\Ri^{n+m}$ all $t>0$ and all $y\in B(x\,;t^{1/2})$.
Thus (\ref{epi4.41}) is valid for all $x,y\in\Ri^{n+m}$ and all $t>0$ with $d(x\,;y)^2/t\leq1$.
Under the latter restraint one can of course introduce a Gaussian factor to obtain the desired lower bound.
Therefore it remains to derive the bound for $d(x\,;y)^2/t\geq1$.
This can be achieved by combination of the semigroup property,  the volume doubling property  and the bound (\ref{epi4.41}) by  adaptation of an argument of Jerison and Sanchez-Call\'e, \cite{JSC} Section~5.

Let $\rho=d(x\,;y)$ and  assume $\rho^2/t\geq 1$. 
Choose a  continuous path of length $l\leq 2\,\rho$ connecting $x$ and $y$.
Next let $k\geq 4$ be the  integer   satisfying $k\geq 4\rho^2/t> k-1$.
Then fix points $x_1,\ldots,x_{k-1}$ in the path with 
$d(x_j\,;x_{j+1})\leq 2\rho/k$ for $j\in\{0,\ldots, k-1\}$ where $x_0=x$ and $x_k=y$.
Now set $B_j=B(x_j\,;2\rho/k)$ and $I_j=B_j\cap B_{j+1}$.
Then
\begin{equation}
K_t(x\,;y)\geq \int_{I_1\times\ldots\times I_{k-1}} \hspace{-4mm}d\xi_1\ldots d\xi_{k-1}\,K_{t/k}(x\,;\xi_1)\,K_{t/k}(\xi_1\,;\xi_2)\ldots K_{t/k}(\xi_{k-1}\,;y)
\;.
\label{epi4.411}
\end{equation}
If $\xi_j\in I_j$ and $\xi_{j+1}\in I_{j+1}$ then $d(\xi_j\,;\xi_{j+1})\leq 2\rho/k$.
Thus since $4\rho^2/t\leq k$  one has $d(\xi_j\,;\xi_{j+1})^2/(t/k)\leq (4\rho^2/t)/k\leq1$.
Therefore it follows from (\ref{epi4.41}) that 
\[
K_{t/k}(\xi_j\,;\xi_{j+1})\geq (c/a)|B(\xi_j\,;(t/k)^{1/2})|^{-1}\;.
\]
Moreover, for each $\xi_j\in B(x_j\,; (t/k)^{1/2})$ one has $B(\xi_j\,;(t/k)^{1/2})\subseteq B(x_j\,;2(t/k)^{1/2})$.
Therefore
\[
|B(\xi_j\,;(t/k)^{1/2})|\leq |B(x_j\,;2(t/k)^{1/2})|\leq b\,|B(x_j\,;(t/k)^{1/2})|
\]
 where the second step uses the volume doubling property (\ref{epi4.2}).
Hence 
\[
K_{t/k}(\xi_j\,;\xi_{j+1})\geq  (c/ab)|B(x_j\,;(t/k)^{1/2})|^{-1}
\]
for each $j\in\{0,\ldots,k-1\}$.
Further $I_j$ contains a ball $B(\xi\,;\rho/k)$ with $B(x_{j}\,;2\rho/k)\subseteq B(\xi\,;4\rho/k)$.
Since $|B(\xi\,;4\rho/k)|\leq b^{2}\,|B(\xi\,;\rho/k)|$ by (\ref{epi4.2}) one then has
\[
|I_j|\geq   b^{-2}\,|B(x_{j}\,;2\rho/k)|\geq b^{-2}\,|B(x_{j}\,;(t/4k)^{1/2})|\geq b^{-3}|B(x_{j}\,;(t/k)^{1/2})|
\;.
\]
The second  inequality follows because $4\rho^2/t> k-1$.
Hence $(2\rho/k)^2\geq t/4k$.
The third uses volume doubling.
Combination of these estimates then gives
\begin{eqnarray*}
K_t(x\,;y)&\geq &(c/ab)^k\Big(\prod^{k-1}_{j=0}|B(x_j\,;(t/k)^{1/2})|^{-1}\Big)
\,b^{-3(k-1)}\,\Big(\prod^{k-1}_{j=1}|B(x_j\,;(t/k)^{1/2})|\Big)\\[5pt]
&=& (c/ab)(c/ab^4)^{k-1}|B(x\,;(t/k)^{1/2})|^{-1}\geq (c/ab)(c/ab^4)^{k-1}|B(x\,;t^{1/2})|^{-1}
\;.
\end{eqnarray*}
Since $k-1< 4\rho^2/t$ one then obtains the lower bounds
\[
K_t(x\,;y)\geq (c/ab)\,|B(x\,;t^{1/2})|^{-1}\,e^{-\omega d(x;y)^2/t}
\]
with $\omega=\log(ab^4/c)$
for all $x,y\in \Ri^{n+m}$ and $t>0$ with $d(x;y)^2/t\geq1$ .
This completes the proof of Theorem~\ref{tpi1.2}.
\hfill$\Box$

\bigskip

\noindent{\bf Proof of Theorem~\ref{tpi1.3}}$\;$The proof is very similar but it relies on Sturm's extension \cite{Stu4, Stu5}  of 
Grigor'yan and Saloff-Coste's work characterizing   the parabolic Harnack inequality.
Sturm establishes that the parabolic Harnack inequality holds for a large class of strictly local regular Dirichlet spaces $X$ with an intrinsic distance $\rho$ if the volume doubling property and the Poincar\'e inequality are satisfied.
The key point is that the intrinsic distance is finite, continuous, defines the original topology of the space and $(X,\rho)$ is complete. 
(For an extensive  discussion in a setting similar to ours see \cite{GyS} and especially Theorem~2.31.)
Therefore the proof of Theorem~\ref{tpi1.3}  reduces to verifying the assumptions of Sturm's theorem for the Dirichlet forms $h_\pm$ forms associated with the generators $H_\pm$ of the submarkovian semigroups $S^{(\pm)}$
on the spaces $L_2(\Omega_\pm)$ and for  the distance functions $\rho_\pm$ obtained by restricting the Riemannian
distance $d(\,\cdot\,;\,\cdot\,)$ to $\Omega_\pm$.

The Dirichlet forms $h_\pm$ are, however, clearly strictly local and regular. 
Thus it remains to consider  properties of the Riemannian distance $d(\,\cdot\,;\,\cdot\,)$.
Let $d_e(\,\cdot\,;\,\cdot\,)$ denote the standard Euclidean distance on $\Ri^{n+m}$.
Since  the Riemannian distance is equivalent to the distance $D_\delta(\,\cdot\,;\,\cdot\,)$  given by  (\ref{epi3.41})
it follows  that for each  $x\in \Ri^{n+m}$ and $r>0$ there exists
a positive real number $r'>0$ such that $d_e(x\,;y) < r'$ implies that $d(x\,;y)< r$
and conversely $d(x\,;y) <r$ implies $d_e(x\,;y)< r'$.
Thus  the distances $d_e(\,\cdot\,;\,\cdot\,)$ and $d(\,\cdot\,;\,\cdot\,) $ 
determine the same topology.
Then a  standard argument establishes that 
$ \Ri^{n+m}$ and any of its closed subsets are complete with respect to both distances. 
It follows  that the spaces $ \Ri^{n+m}$ and $ \Ri_\pm\times \Ri^{m} \subset  \Ri^{m+1}$ satisfy assumptions (A1) and 
(A2)  of \cite{GyS}, page 24. Note that this implies that these spaces are geodesic length spaces in the terminology
of Theorem 2.11 of \cite{GyS}.

Finally these observations establish that  the theorem of Sturm applies to the operators $H_\pm$.
Hence they  satisfy the parabolic Harnack inequality on $\Omega_\pm$.
Then the proof of Theorem~\ref{tpi1.3} is  a repetition of the arguments used to establish 
Theorem~\ref{tpi1.2}. \hfill$\Box$

\bigskip

Theorems~\ref{tpi1.2} and \ref{tpi1.3} demonstrate that the heat semigroup corresponding to the degenerate operator $H$ has
a Gaussian character similar to that of a non-degenerate strongly elliptic operator.
Even in the non-ergodic situation $n=1$, $\delta_1\in[1/2,1\rangle$, the Gaussian characteristics persist in the ergodic components.
The Gaussian upper bounds on the kernel are, however, not  optimal.
These estimates can be improved as in the strongly elliptic case, e.g.\ for each $\varepsilon>0$ one can choose $a'$ such that $\omega'=(4+\varepsilon)^{-1}$.

\section{The exceptional case}\label{S5}

The discussion of the heat kernel in Section~\ref{S4} does not cover the  case $n=1$, 
$\delta_1\in[0,1/2\rangle$ and $\delta_1'\in[1/2,1\rangle$.
Moreover, it follows from Theorem~\ref{tpi1.1}.II that  in this case the uniform Poincar\'e inequality is not valid.
The proof in Section~\ref{S3} that the inequality is invalid demonstrates that the problem is a global one.
In this section we establish that the inequality is  nevertheless valid locally
and subsequently discuss the implications for the heat semigroup.

The principal result is the following.
\begin{thm}\label{tpi5.1}
Assume  $n=1$,  $\delta_1\in[0,1/2\rangle$ and $\delta_1'\in[1/2,1\rangle$.
Then there is a $\kappa\in\langle0,1]$ and for each $R>0$ there is a $\lambda_R>0$  such that 
\begin{equation}
\int_{B(x;r)}dy\,\Gamma(\varphi)(y)\geq \lambda_R\,r^{-2}\int_{B(x;\kappa r)}dy\left(\varphi(y)-\langle\varphi\rangle_{B}\right)^2
\label{epi5.5}
\end{equation}
for all $x\in \Ri^{1+m}$, $r\in\langle0,R]$ and  $\varphi\in C^1(\Ri^{n+m})$ where $\langle\varphi\rangle_{B}=|B(x\,;\kappa r)|^{-1}\int_{B(x;\kappa r)}dy\,\varphi(y)$.
\end{thm}

The conclusion of the theorem is considerably weaker than Statement~I of Theorem~\ref{tpi1.1} since $\lambda_R$  tends to zero
as $R\to\infty$.
In fact the proof of the theorem establishes that the rate of convergence is given by a power of $R^{-\alpha'(\delta_1'-\delta_1)}$.

The proof is similar to the proof (\ref{epi1.5}).
It consists of three steps.
First one proves the  inequality (\ref{epi5.5}) for balls centred at the origin,
secondly  for balls which do not contain the origin and finally one deduces the result for general balls
from the two special cases.
The only major change occurs in the first step, the discussion of balls centred at the origin.
The essential feature is the following analogue of Proposition~\ref{ppi3.2}.

\begin{prop}\label{ppi5.2}
Assume  $n=1$,  $\delta_1\in[0,1/2\rangle$ and $\delta_1'\in[1/2,1\rangle$.
Then for each $T>0$ there is a $\lambda_T>0$ such that 
\begin{equation}
\int_{C_t}dx\,\Gamma_\delta(\varphi)(x)\geq\lambda_T\,t^{-2}\int_{C_t}dx\,(\varphi(x)-\langle\varphi\rangle)^2
\label{epi4.22}
\end{equation}
for all $\varphi\in C_c^1(\Ri^{1+m})$  and $t\in\langle0,T\,]$ where $\langle\varphi\rangle=|C_t|^{-1}\int_{C_t}dx\,\varphi(x)$.
\end{prop}
\proof\
First if $t\leq 1$ then the integral on the left hand side of (\ref{epi4.22}) is independent of $\delta_1'$ and  the inequality is 
a corollary of Proposition~\ref{ppi3.2}.

Secondly if $t\geq 1$ then $C_t$ is the product of an interval $-t^{\alpha'}\leq x_1\leq t^{\alpha'}$ and a cube
$|x_2|_\infty\leq t^{\beta'}$ in $\Ri^m$.
Then changing variables to $y_1=t^{-\alpha'}x_1$ and $y_2=t^{-\beta'}x_2$  and setting $\psi(y)=\varphi(x)$ one has
\begin{eqnarray*}
\int_{C_t}dx\,\Gamma_\delta(\varphi)(x)&=&
t^{\alpha'+m\beta'}\int_{C_1}dy\,\Big( t^{-2\alpha'}\,|t^{\alpha'}y_1|^{(2\delta_1,2\delta_1')}\,
|(\partial_{y_1}\psi)(y)|^2\\[5pt] &&\hspace{5cm}{}+  t^{-2\beta'}\,|t^{\alpha'}y_1|^{(2\delta_2,2\delta_2')}\,
|(\nabla_{y_2}\psi)(y)|^2\Big)\;.
\end{eqnarray*}
Since $\delta_1'>\delta_1$, $t\geq1 $ and $|y_1|\leq 1$ it follows from Proposition~\ref{ppi2.1} that
\[
 t^{-2\alpha'}\,|t^{\alpha'}y_1|^{(2\delta_1,2\delta_1')}\geq 2^{-4\delta_1'}\, t^{-2\alpha'}\,t^{2\alpha'\delta_1}\,|y_1|^{2\delta_1}
\geq 2^{-4}\,t^{-2}\,T^{-2\alpha'(\delta_1'-\delta_1)}
\]
for all $t\in\langle0,T]$ where we have used $1-\alpha'+\alpha'\delta_1=-\alpha'(\delta_1'-\delta_1)<0$.
Similarly
\[
t^{-2\beta'}\,|t^{\alpha'}y_1|^{(2\delta_2,2\delta_2')}\geq 2^{-2(\delta_2+\delta_2')}\,t^{-2\beta'}\,t^{2\alpha'(\delta_2\vee\delta_2')}\,|y_1|^{2\delta_2}
\geq 2^{-2(\delta_2+\delta_2')}\,t^{-2}\,|y_1|^{2\delta_2}
\]
where the last step uses $1-\beta'+\alpha'(\delta_2\vee\delta_2')=\alpha'(\delta_2\vee\delta_2'-\delta_2')\geq0$.
Combining these estimates one has
\begin{eqnarray*}
\int_{C_t}dx\,\Gamma_\delta(\varphi)(x)&\geq& a_\delta(T)\,t^{-2}\,t^{\alpha'+m\beta'}\int_{C_1}dy\,
\Big(|y_1|^{2\delta_1}\,
|(\partial_{y_1}\psi)(y)|^2+ |y_1|^{2\delta_2}\,
|(\nabla_{y_2}\psi)(y)|^2\Big)\\[5pt]
&=&a_\delta(T)\,t^{-2}\,t^{\alpha'+m\beta'} \int_{C_1}dy\,\Gamma_\delta(\psi)(y)
\end{eqnarray*}
for all $t\in\langle0,T]$ where $a_\delta(T)=2^{-4}\,T^{-2\alpha'(\delta_1'-\delta_1)}\,(1\wedge 2^{-2(\delta_2+\delta_2'-2)})$.
But $C_1=[-1,1]^{1+m}$ so it follows from Proposition~\ref{ppi3.2} that there is a $\lambda>0$ such that 
\begin{eqnarray*}
t^{\alpha'+m\beta'}\int_{C_1}dy\,\Gamma_\delta(\psi)(y)&\geq& \lambda\,t^{\alpha'+m\beta'}\int_{C_1}dy\,(\psi(y)-\langle\psi\rangle)^2\\[5pt]
&=&\lambda \int_{C_t}dx\,(\varphi(x)-\langle\varphi\rangle)^2
\end{eqnarray*}
where the last step  follows by reverting  to the original $x$-coordinates.
Finally one concludes by combination of these estimates that (\ref{epi4.22}) is valid with $\lambda_T =\lambda\,a_\delta(T)$.
\hfill$\Box$

\bigskip

Now the proof of Theorem~\ref{tpi5.1} is essentially a corollary of Proposition~\ref{ppi5.2} and the arguments used to prove Theorem~\ref{tpi1.1}.

\medskip

\noindent{\bf Proof of Theorem~\ref{tpi5.1}}$\;$
First, it follows from Proposition~\ref{ppi5.2} and and the  embedding statements, Lemmas~\ref{lpi3.3} and \ref{lpi3.4}, that there is a $\kappa\in\langle0,1]$ and for each $R>0$ there is a $\lambda_R>0$ such that 
the Poincar\'e inequality (\ref{epi5.5})  is valid for all balls $B_\Delta(0\,;r)$ with $r\in\langle0,R\,]$.
Note that the embedding lemmas are general geometric results which are valid for all $n$  all $\delta_1,\delta_1'\in[0,1\rangle$
and balls of arbitrary radius.
Moreover, the value of $\kappa$ which occurs in Lemma~\ref{lpi3.4} is independent of  the radius of the balls.
It also follows from these  lemmas together with the estimates in the foregoing proof that $\lambda_R$ converges to zero as $R\to\infty$ and the rate of convergence is given by a power of $R^{-\alpha'(\delta_1'-\delta_1)}$.

Secondly,  it follows from the discussion of Case~II of the proof of Theorem~\ref{tpi1.1}.I that one may also choose
$\kappa$ and $\lambda_R$ such that  (\ref{epi5.5}) is valid for all balls $B_\Delta((\xi_1,0)\,;r)$ with $r<r_\xi$.
Note that the arguments in Case~II are independent of the assumption $\delta_1\vee\delta_1'\in[0,1/2\rangle$ if $n=1$
(see Remark~\ref{rpi3.1}).
Hence the arguments apply with  the current assumptions $\delta_1\in[0,1/2\rangle$ and $\delta_1'\in[1/2, 1\rangle$.

It now remains to establish (\ref{epi5.5}) for the balls $B_\Delta((\xi_1,0)\,;r)$ with $r_\xi<R$ and $r\in [r_\xi, R\,]$.
But this is again a corollary of the foregoing special cases.

Set $K=2\,(1+\kappa)/\kappa$.
First suppose $K\,r_\xi\leq R$.
Then if $r\in [K\,r_\xi,R\,]$ the Poincar\'e inequality follows from the previous two special cases by the first argument in the discussion of Case~III of the proof of Theorem~\ref{tpi1.1}.I.
If, however, $r\in[r_\xi, Kr_\xi]$ it follows from the special cases
by the second argument in the discussion of Case~III.
Secondly, if  $R\leq K\,r_\xi$ then the condition $r\in[r_\xi,R\,]$ implies that $r\in[r_\xi, K\,r_\xi]$ and the Poincar\'e inequality follows again.

Finally combination of these conclusions gives the Poincar\'e inequality (\ref{epi5.5}) as stated in Theorem~\ref{tpi5.1}.
\hfill$\Box$

\bigskip

The Poincar\'e inequality (\ref{epi5.5}) of Theorem~\ref{tpi5.1} is a local version of the inequality (\ref{epi1.5})
established in Theorem~\ref{tpi1.1} for $n\geq2$ or for $n=1$ and $\delta_1\vee\delta_1'\in[0,1/2\rangle$
insofar the value of $\lambda$ depends on the scale $R$ of the balls.
But the results of Grigor'yan and Saloff-Coste establish that much of the discussion of Section~\ref{S4} still applies 
to the heat kernel although the conclusions are of a local nature.
In particular  the volume doubling property in combination with the 
 local Poincar\'e inequality gives a local version of the parabolic Harnack inequality.
Explicitly, if $n=1$,  $\delta_1\in[0,1/2\rangle$ and $\delta_1'\in[1/2,1\rangle$ then for each $R>0$ 
there exists an $a>0$
 such that for any $x\in\Ri^{1+m}$,  $t>0$ and all $r\in\langle0,R\,]$ any non-negative (weak) solution $\varphi$ of the parabolic equation $(\partial_t+H)\varphi=0$ in the cylinder $Q=\langle t,\, t+r^2\rangle\times B(x\,;2r)$
satisfies
\begin{equation}
\sup_{Q^-}\varphi\leq a\,\inf_{Q^+}\varphi
\label{epi5.1}
\end{equation}
where $Q^-=[\,t+r^2\!/4,\, t+r^2\!/2\,]\times B(x\,;r)$ and $Q^+=[\, t+3\hspace{1pt} r^2\!/4, \,t+r^2\rangle\times B(x\,;r)$.

The local version of the Harnack inequality again implies H\"older continuity of the heat kernel by Moser's argument.
Therefore Gaussian upper bounds on the kernel follow  from the almost everywhere bounds of Corollary~6.6 in \cite{RSi2}.

\begin{cor}\label{cpi5.1}$\;$Assume  $n=1$, 
$\delta_1\in[0,1/2\rangle$ and $\delta_1'\in[1/2,1\rangle$.
Then the semigroup kernel $x,y\in\Ri^{1+m}\mapsto K_t(x\,;y)$ is jointly H\"older continuous 
and there exist $a',\omega'>0$ such that 
\[
 K_t(x\,;y)\leq a'\,|B(x\,;t^{1/2})|^{-1}\,e^{-\omega'\,d(x;y)^2/t}
\]
 for all $x,y\in\Ri^{1+m}$ and $t>0$.
\end{cor}

Although the heat kernel satisfies  Gaussian upper bounds this type of bound is not expected to be optimal for reasons we discuss below.
In addition one cannot  expect matching Gaussian lower bounds as these would imply the global parabolic Harnack inequality
which in turn would imply the global Poincar\'e inequality in contradiction with Theorem~\ref{tpi1.1}.II.
Nevertheless one has the on-diagonal lower bounds (\ref{epi4.4}) established in \cite{RSi2} and then
arguing with the local Harnack inequality as in the proof of Theorem~\ref{tpi1.2} one obtains the small time off-diagonal lower bound
\begin{equation}
K_t(x\,;y)\geq (c/a)\,|B(x\,;t^{1/2})|^{-1}
\label{epi5.41}
\end{equation}
with  the restrictions  $d(x\,;y)\leq t^{1/2}\leq R$
where $a$ and $R$ are the parameters in~(\ref{epi5.1}).
These small-time lower bounds then  imply that the kernel $K_t$ is strictly positive for all $t>0$.
This is a consequence of the semigroup property and an estimate of the type (\ref{epi4.411}).
Consequently the semigroup $S$ is ergodic, i.e.\ there are no non-trivial $S$-invariant subspaces of the form $L_2(\Omega)$.

The complications with this exceptional case arise because the subspaces $\Omega_\pm$ are `approximately' invariant.
We will not discuss the precise meaning of approximately invariant but 
instead argue that there is  a similarity of the evolution in the exceptional case and the evolution on
 manifolds with ends  as described by Grigor'yan and Saloff-Coste  \cite{GSC}.
In the special case of manifolds with two   ends  with the same dimension the situation 
 can be described as two copies of $\Ri^n$ connected by a compact cylinder. 
Assuming the manifold is rotationally  invariant it can be identified as
 $\Ri \times \Si^n$ with polar coordinates $(r, \sigma) \in  \Ri \times \Si^n$
the Riemannian metric is given by $ d_r^2 + f(r) d_\sigma^2$, where 
$f(r)>0$ is continuous and $f(r)=r^{-2}$ for $|r| \ge 1$. Then the  quadratic form 
corresponding to the Laplace Beltrami operator can be defined as 
\[
Q(\psi) =\int_{\Ri}\int_{\Si^n} (|\partial_r\psi|^2 +f(r)^{-1}  |\partial_\sigma\psi|^2)g(r)drd\sigma 
\]
where $g(r)>0 $ is continuous and $g(r)=r^{n-1}$ for $|r| \ge 1$.
If one is just interested  in the evolution corresponding to Brownian motion in the radial direction 
then we can restrict attention to functions which are invariant in~$\sigma$. 
This 
leads to the  one-dimensional operator $L$ acting on $L_2(\Ri\,;g(r)dr)$  corresponding 
to the quadratic form 
\[
Q_1(\phi)=\int_{\Ri}|\partial_r\psi|^2g(r)dr
\;.
\]
After  a simple change  of variable, which we describe in Example~\ref{expi5.1} below, the operator $L$ is equivalent to a 
one-dimensional degenerate elliptic operator. 
It is worth noting that if we consider  an operator $H$ acting on $L^2(\Ri\times \Ri^m)$
and are only interested  in the evolution of $x_1$ then we  again obtain the same operator
or rather its equivalent version discussed in Example~\ref{expi5.1}. 
This means that at least in the radial case the approximate invariance corresponds is characterized by the  heat kernel bounds described in  \cite{GSC}.

We  conclude with an example, adapted from \cite{HSi},  which illustrates the structure of the 
$H_{\delta,0}$ and $S^{\delta,0}$ in the simplest case $\delta_1=0$
and the connection with the end problem.

\begin{exam}\label{expi5.1}
Let $H_0=-d_x\,(1\vee|x|)^{2\delta}\,d_x$ with $\delta\in[1/2,1\rangle$ be the operator on $L_2(\Ri)$ with domain $C_c^\infty(\Ri)$.
Since the coefficient of $H_0$ is strictly positive the operator is essentially self-adjoint.
Let $H$ denote the self-adjoint closure and $h$ the corresponding Dirichlet form.
Now for $\varphi\in C_c^\infty(\Ri)$ define $\Phi=\varphi\circ f$ where $f(x)=(-|x|^\alpha) \wedge \,x$ if $x\leq 0$ and 
$f(x)=\,x\vee|x|^\alpha$ if $x>0$ 
with $\alpha=(1-\delta)^{-1}\in[2,\infty\rangle$.
The mapping $\varphi\mapsto \Phi$ extends to an isometric isomorphism $U$ from $L_2(\Ri)$ to $L_2(\Ri\,;\mu)$ where
$d\mu=f'$.
Thus the $L_2(\Ri\,;\mu)$-norm $\|\cdot\|_{2,\mu}$ is given by
\begin{eqnarray*}
\|\Phi\|_{2,\mu}^2&=&\int_\Ri dy\,f'(y)\,|\Phi(y)|^2\\[5pt]
&=&\alpha\int_{-\infty}^{-1}dy\,|y|^{\alpha-1}\,|\Phi(y)|^2
+\int^1_{-1}dy\,|\Phi(y)|^2+
\alpha\int_1^\infty dy\,|y|^{\alpha-1}\,|\Phi(y)|^2
\;.
\end{eqnarray*}
But if $\alpha=k$ is a positive integer, i.e.\ if $\delta=1-k^{-1}$, then $dy\, |y|^{k-1}$ is the radial measure on $\Ri^k$.
Therefore the first and 
 last integrals can be identified as the square of the radial part of the $L_2(\Ri^k)$-norm   
 restricted to $\Ri^k\backslash B(0\,;1)$ where $B(0\,;1)$ denotes the unit Euclidean ball centred at the origin.
Next one calculates that $h(\varphi)=h_\mu(\Phi)$ for each $\varphi\in C_c^\infty(\Ri)$ where 
\begin{eqnarray*}
h_\mu(\Phi)
&=&\alpha^{-1}\int_{-\infty}^{-1}dy\,|y|^{\alpha-1}\,|\Phi'(y)|^2
+\int^1_{-1}dy\,|\Phi'(y)|^2+
\alpha^{-1}\int_1^\infty dy\,|y|^{\alpha-1}\,|\Phi'(y)|^2
\;.
\end{eqnarray*}
This identification then extends by closure to all $\varphi\in D(h)$.
The form $h_\mu$ is a Dirichlet form on $L_2(\Ri\,;\mu)$ with $D(h_\mu)=U D(h)$.
Therefore  the corresponding self-adjoint operator $H_\mu= UHU^{-1}$ generates the  submarkovian semigroup 
 $S^\mu_t=US_t U^{-1}$.
If  $\alpha=k$ is a positive integer the operator $H_\mu$ models the radial part of the Laplace-Beltrami operator acting on 
a manifold which consists of two copies  $\Ri^k\backslash B(0\,;1)$ with  a cylindrical channel joining the two balls.
The   semigroup $S^\mu$  describes a diffusion process for which the two ends  $\Ri^k\backslash B(0\,;1)$ are largely invariant since the probability of passing from one end of the manifold to the other is small.
Such processes have been studied  by Grigor'yan and Saloff-Coste \cite{GSC}.
In particular they have derived matching upper and lower bounds which describe  non-Gaussian behaviour.
\end{exam}

\bigskip

\subsection*{Acknowledgement}
This collaboration  was carried out during numerous visits of the second author to the Mathematical Sciences Institute at ANU with the support of ARC Discovery grant DP130101302.


\begin{thebibliography}{FGW94}

\bibitem[BCF96]{BCF}
{\sc Benjamini, I., Chavel, I., {\rm and} Feldman, E.~A.}, Heat kernel lower
  bounds on {R}iemannian manifolds using the old ideas of {N}ash.
\newblock {\em Proc.\ London Math.\ Soc.} {\bf 72} (1996),  215--240.

\bibitem[CF91]{ChF}
{\sc Chavel, I., {\rm and} Feldman, E.~A.}, Isoperimetric constants, the
  geometry of ends, and large time heat diffusion in {R}iemannian manifolds.
\newblock {\em Proc.\ London Math.\ Soc.} {\bf 62} (1991),  427--448.

\bibitem[Dav97]{Dav15}
{\sc Davies, E.~B.}, Non-{G}aussian aspects of heat kernel behaviour.
\newblock {\em J.\ London Math.\ Soc.} {\bf 55} (1997),  105--125.

\bibitem[FGW94]{FGW2}
{\sc Franchi, B., Guti{\'e}rrez, C.~E., {\rm and} Wheeden, R.~L.}, Weighted
  Sobolev--Poincar\'e inequalities for Grushin type operators.
\newblock {\em Comm.\ Part.\ Diff.\ Eq.} {\bf 19} (1994),  523--604.

\bibitem[FKS82]{FKS}
{\sc Fabes, E.~B., Kenig, C.~E., {\rm and} Serapioni, R.~P.}, The local
  regularity of solutions of degenerate elliptic equations.
\newblock {\em Comm.\ Part.\ Diff.\ Eq.} {\bf 7} (1982),  77--116.

\bibitem[Gri92]{Gri4}
{\sc Grigor'yan, A.}, The heat equation on noncompact Riemannian manifolds.
\newblock {\em Math.\ USSR-Sb.} {\bf 72}, No.\ 1 (1992),  47--77.
\newblock Mat.\ Sb.\ {\bf 182}, No.\ 1, (1991), 55--87.

\bibitem[Gru70]{Gru}
{\sc Gru\v{s}in, V.~V.}, A certain class of hypoelliptic operators.
\newblock {\em Mat.\ Sb.\ (N.S.)} {\bf 83 $($125$)$} (1970).

\bibitem[GSC09]{GSC}
{\sc Grigor'yan, A., {\rm and} Saloff-Coste, L.}, Heat kernel on manifolds with
  ends.
\newblock {\em Ann.\ Inst.\ Fourier, Grenoble} {\bf 59} (2009),  1917--19974.

\bibitem[GSC11]{GyS}
{\sc Gyrya, P., {\rm and} Saloff-Coste, L.}, {\em Neumann and {D}irichlet heat
  kernels in inner uniform domains}, vol.\ 336.
\newblock Soci\'et\'e Math\'ematique de France, Paris, 2011.

\bibitem[H{\"o}r67]{Hor1}
{\sc H{\"o}rmander, L.}, Hypoelliptic second order differential equations.
\newblock {\em Acta Math.} {\bf 119} (1967),  147--171.

\bibitem[HS09]{HSi}
{\sc Hassell, A., {\rm and} Sikora, A.}, Riesz transforms in one-dimension.
\newblock {\em Ind.\ Univ.\ Math.\ J.} {\bf 58} (2009),  823--852.

\bibitem[Jer86]{Jer}
{\sc Jerison, D.}, The Poincar{\'e} inequality for vector fields satisfying
  H{\"o}rmander's condition.
\newblock {\em Duke Math.\ J.} {\bf 53} (1986),  503--523.

\bibitem[JSC86]{JSC}
{\sc Jerison, D.~S., {\rm and} S{\'a}nchez-Calle, A.}, Estimates for the heat
  kernel for a sum of squares of vector fields.
\newblock {\em Ind.\ Univ.\ Math.\ J.} {\bf 35} (1986),  835--854.

\bibitem[Lu94]{Lu2}
{\sc Lu, G.}, The sharp {P}oincar\'e inequality for free vector fields: an
  endpoint result.
\newblock {\em Revista Matem\'atica Iberoamericana} {\bf 10} (1994),  453--466.

\bibitem[Mos61]{Mos0}
{\sc Moser, J.}, On {H}arnack's theorem for elliptic differential equations.
\newblock {\em Commun.\ Pure Appl.\ Math.} {\bf 14} (1961),  577--591.

\bibitem[Mos64]{Mos}
\leavevmode\vrule height 2pt depth -1.6pt width 23pt, A Harnack inequality for
  parabolic differential equations.
\newblock {\em Commun.\ Pure Appl.\ Math.} {\bf 17} (1964),  101--134.
\newblock Correction to: ``A Harnack inequality for parabolic differential
  equations'', Comm.\ Pure Appl.\ Math.\ {\bf 20} (1967), 231--236.

\bibitem[Mos71]{Mos1}
\leavevmode\vrule height 2pt depth -1.6pt width 23pt, On pointwise estimate for
  parabolic differential equations.
\newblock {\em Commun.\ Pure Appl.\ Math.} {\bf 24} (1971),  727--740.

\bibitem[MR92]{MR}
{\sc Ma, Z.~M., {\rm and} R{\"o}ckner, M.}, {\em Introduction to the theory of
  $($non symmetric$)$ Dirichlet Forms}.
\newblock Universitext. Springer-Verlag, Berlin etc., 1992.

\bibitem[RS08]{RSi2}
{\sc Robinson, D.~W., {\rm and} Sikora, A.}, Analysis of degenerate elliptic
  operators of Gru\v{s}in type.
\newblock {\em Math.\ Z.} {\bf 260} (2008),  475--508.

\bibitem[SC92a]{Sal2}
{\sc Saloff-Coste, L.}, A note on Poincar{\'e}, Sobolev, and Harnack
  inequalities.
\newblock {\em Internat.\ Math.\ Res.\ Notices} {\bf 1992}, No.\ 2 (1992),
  27--38.

\bibitem[SC92b]{Sal3}
\leavevmode\vrule height 2pt depth -1.6pt width 23pt, Uniformly elliptic
  operators on Riemannian manifolds.
\newblock {\em J. Diff.\ Geom.} {\bf 36} (1992),  417--450.

\bibitem[SC95]{Sal4}
\leavevmode\vrule height 2pt depth -1.6pt width 23pt, Parabolic Harnack
  inequality for divergence-form second-order differential operators.
\newblock {\em Potential Anal.} {\bf 4} (1995),  429--467.

\bibitem[Stu95]{Stu4}
{\sc Sturm, K.-T.}, Analysis on local Dirichlet spaces. II. Upper Gaussian
  estimates for the fundamental solutions of parabolic equations.
\newblock {\em Osaka J. Math.} {\bf 32} (1995),  275--312.

\bibitem[Stu96]{Stu5}
\leavevmode\vrule height 2pt depth -1.6pt width 23pt, Analysis on local
  Dirichlet spaces. III. The parabolic Harnack inequality.
\newblock {\em J. Math.\ Pures Appl.} {\bf 75} (1996),  273--297.

\bibitem[Tru73]{Tru2}
{\sc Trudinger, N.~S.}, Linear elliptic operators with measurable coefficients.
\newblock {\em Ann.\ Scuola Norm.\ Sup.\ Pisa} {\bf 27} (1973),  265--308.

\end{thebibliography}
\end{document}